\documentclass[aap,MSNbibl,seceqn,dvips]{arximspdf}
\usepackage{accents}

\doi{10.1214/12-AAP874} 
\volume{23}
\issue{3}
\pubyear{2013}
\firstpage{1254}
\lastpage{1289}

\makeatletter
\newcommand{\rrvert}{\vert}
\newcommand{\llvert}{\vert}
\renewcommand{\mathring}[1]{\accentset{\circ}{#1}}
\newcommand{\E}{\mathbf{E}}
\newcommand{\PP}{\mathbf{P}}

\newcommand{\Bin}{\operatorname{Bin}}

\newtheorem{tm}{Theorem}[section]

\newtheorem{lem}{Lemma}
\newtheorem{col}{Corollary}
\newproclaim{remark}{Remark}
\newproclaim{example}{Example}
\makeatother

\begin{document}
\begin{frontmatter}

\title{Degree and clustering coefficient in sparse random intersection graphs}
\runtitle{Degree and clustering}

\begin{aug}
\author{\fnms{Mindaugas}~\snm{Bloznelis}\corref{}\thanksref{t1}\ead[label=e1]{mindaugas.bloznelis@mif.vu.lt}\ead[label=u1,url]{http://www.mif.vu.lt/\textasciitilde bloznelis}}
\thankstext{t1}{Supported by Research Council of Lithuania Grant MIP-053/2011.}
\runauthor{M. Bloznelis}
\affiliation{Vilnius University}
\dedicated{Dedicated to Professor Friedrich G\"otze on the occasion of his 60th birthday}
\address{Faculty of Mathematics and Informatics\\
Vilnius University\\
Naugarduko 24\\
Vilnius LT-03225\\
Lithuania\\
\printead{e1}\\
\printead{u1}}
\end{aug}

\received{\smonth{10} \syear{2011}}
\revised{\smonth{5} \syear{2012}}

%
\begin{abstract}
We establish
asymptotic vertex degree distribution and examine its relation to the
clustering coefficient
in two popular
random intersection graph models of Godehardt and Jaworski [\textit
{Electron. Notes Discrete Math.} \textbf{10} (2001) 129--132].
For sparse graphs with a positive clustering coefficient, we examine
statistical dependence between the (local) clustering coefficient
and the degree.
Our results are mathematically rigorous.
They
are consistent with the empirical observation of Foudalis et al.
[In \textit{Algorithms and Models for Web Graph} (2011) Springer]
that, ``clustering correlates negatively with degree.''
Moreover, they explain empirical results on $k^{-1}$ scaling of the
local clustering coefficient of a vertex of degree $k$ reported
in Ravasz and Barab\'asi [\textit{Phys. Rev. E} \textbf{67} (2003) 026112].
\end{abstract}

%

\end{frontmatter}

\section{Introduction}\label{sec1}
In a recent paper~\cite{Foudalis2011}, Foudalis et al. analyzed
Facebook data and made the observation
that ``clustering correlates negatively with degree.'' Their empirical findings
suggest that the chances of two neighbors of a given vertex to be
adjacent is a
decreasing function of the vertex degree.
A reasonable question to ask is whether and how such a phenomenon
can be explained with the aid of a known
theoretical
model.

This question is addressed in the present paper.
We consider two simple random
graph models admitting a power-law degree distribution and positive
clustering coefficient.

Given a finite set $W$ and a collection
of its subsets $D_1,\ldots, D_n$,
the active intersection graph defines adjacency relation between the subsets
by declaring\vadjust{\goodbreak} two subsets adjacent whenever they share at least $s$
common elements.
The passive intersection graph defines adjacency
relation between elements of $W$.
A pair of elements is declared an edge if it is contained in $s$ or
more subsets.
Here $s\ge1$ is a model parameter.
Both models have reasonable interpretations: in the active graph, two
actors $D_i$ and $D_j$
establish a communication link whenever they have sufficiently many
common interests;
in the passive graph, students become acquaintances if they participate
in sufficiently
many joint projects. In order to model active and passive graphs with desired
statistical properties, we choose subsets $D_1,\ldots, D_n$ at random
and obtain
random intersection graphs. Such random graph models have been
introduced in
\cite{godehardt2001,karonski1999}; see also~\cite{eschenauer2002}.

We remark
that the adjacency relations of random intersection graphs resemble
that of some
real networks, like, for example, the actor network, where two actors
are linked by an edge whenever they have
acted in the same film, or the collaboration network, where authors are declared
adjacent whenever they have coauthored at least $s$ papers. These
networks exploit the underlying
bipartite graph
structure: actors are linked to films, and authors to papers.
Newman et al.~\cite{Newman+W+S2002} pointed out that the clustering
property of some social networks
could be explained by the presence of
such a bipartite graph structure; see also
\cite{Barbour2011,Guillaume+L2004}, and references therein.
In this respect, it is relevant to mention that
the random intersection graphs
of the present paper
can be obtained from the random bipartite graph with
bipartition $V\cup W$, where
each vertex $v_i$ of the set
$V=\{v_1,\ldots, v_n\}$ selects the set $D_i\subset W$ of its neighbors
in the bipartite graph
independently at random. In addition, we assume for simplicity that
all elements of $W$ have equal probabilities to be selected.
Now, the active intersection graph defines the adjacency relation on
the vertex set $V$:
$v_i$ and $v_j$ are adjacent if
they have at least $s$ common neighbors in the bipartite graph. Similarly,
vertices $w_i,w_j\in W$ are adjacent in the passive graph whenever they
have at least $s$
common neighbors in the bipartite graph.
An attractive property of these models that motivated our study
is that they capture some features of
real networks and, therefore, they may be useful in better
understanding the statistical properties of some social networks.

To summarize our paper, we note that both active and passive random
intersection graphs admit
a nonvanishing clustering coefficient although
their ``clustering mechanisms'' are different.
Both models admit a power-law (asymptotic) degree distribution but of a
different structure.
A common feature of these models is that the (asymptotic) degree
distribution of graphs with a nontrivial clustering
coefficient (i.e., one taking values other than $0$ or $1$) has a
finite second moment.
Another interesting fact is that, in many cases,
the chance of two neighbors of a given vertex to be adjacent decays
as $ck^{-1}$, where $k$ is the vertex degree. For example, the $k^{-1}$
\textit{scaling} is shown for
a power-law active random
intersection graph
with a positive clustering coefficient and integrable degree. This
theoretical result agrees with the
empirical findings
reported in~\cite{RBarabasi2003} showing the $k^{-1}$ scaling in some
real networks.

The paper is organized as follows. In Section~\ref{sec1}, we introduce the
active random intersection
graph and collect results about the degree distribution, clustering
coefficient and
statistical dependence between the degree and the clustering
coefficient of a given vertex. Section~\ref{sec2} is devoted to the passive
random intersection graph.
Proofs are given in Section~\ref{sec3}.

In what follows, the expressions $o(\cdot)$, $O(\cdot)$ refer to the
case where the size $m$ of
the auxiliary set $W$ and the number $n$ of subsets defining the
intersection graph both tend to infinity
(unless stated otherwise).
Given a sequence $\{\xi_m\}$ of random variables, we write $\xi_m=o_P(1)$ if $\xi_m\to0$
in probability, that is, $\forall\varepsilon>0, \PP(\xi_m>\varepsilon
)=o(1)$. We write $\xi_m=O_P(1)$
if the sequence is stochastically bounded, that is,
$\forall\varepsilon>0, \exists A>0$: $\sup_m\PP(|\xi_m|>A)<\varepsilon$.

\section{Active intersection graph}\label{sec2}

Given a set of attributes $W=\{w_1,\ldots, w_m\}$, an actor $v$ is
identified with the set $D(v)$ of
attributes selected by $v$ from $W$. We assume that the actors
$v_1,\ldots, v_n$
choose their attribute sets $D_i=D(v_i)$, $1\le i\le n$,
independently at
random, and we declare $v_i$ and $v_j$ adjacent (denoted $v_i\sim v_j$) whenever
they share at least $s$ common attributes, that is,
$|D_i\cap D_j|\ge s$. Here and below, $s\ge1$ is the same for all
pairs $v_i,v_j$.
The graph on the vertex set $V=\{v_1,\ldots, v_n\}$
defined by this adjacency relation is called the \textit{active}
random intersection graph; see~\cite{godehardt2003}. Subsets of $W$ of
size $s$ play a special role;
we call them joints. They
serve as witnesses of established links: $v_i\sim v_j$ whenever there exists
a joint belonging to both $D_i$ and~$D_j$.

We first assume, for simplicity, that the random sets $D_1,\ldots, D_n$ have
the same probability distribution of the form
%
\begin{equation}
\label{W-D-A} \PP\bigl(|D_i|=A\bigr)=P\bigl(|A|\bigr){\pmatrix{m
\cr
|A|}}^{-1}\qquad
\mbox{for } A\subset W.
\end{equation}
That is, given an integer $k$, all subsets $A\subset W$ of size
$|A|=k$ receive equal chances, proportional to the weight $P(k)$,
where $P$ is a probability on $\{0, 1,\break\ldots, m\}$. The random
intersection graph defined in this way is denoted $G_s(n,\break m,P)$.
One special case where all random sets are of the same (nonrandom)
size
has attracted particular
attention in the literature (see, e.g.,~\cite{Blackburn2009,eschenauer2002,godehardt2003,Rybarczyk2011,Spirakis2011,Yagan2009}),
as it provides a convenient model of a secure wireless network.
We call it the \textit{uniform} active random intersection graph and
denote $G_s(n,m,\delta_x)$.
Here
$\delta_x$ is the
probability distribution putting mass $1$ on $x$, and $x$ is the size
of the random sets.

\subsection{Degree distribution}\label{sec2.1} We consider a sequence
of random intersection graphs, where $m$ and $n=n_m$ tend to infinity,
and where $s=s_m$ and $P=P_m$ depend on $m$.\vadjust{\goodbreak} We suppress the subscript
$m$ in what follows whenever this does not cause an ambiguity. By
$X_i=|D_i|$ we denote the size of the attribute set $D_i$. Note that
$X_i$ is a random variable taking values in $\{0,1,\ldots, m\}$ and
having the probability distribution $P$.

\begin{lem}\label{active-edge} Let $m\to+\infty$. Assume that $s=O(1)$ and
$\E X_1^s{\mathbb I}_{\{X_1>\varepsilon\sqrt{m}\}}=o(\E X_1^s)$ for
any $\varepsilon>0$. Then
for any pair of vertices $v_i,  v_j$, the edge probability in $G_s(n,m,P)$
%
\begin{equation}
\label{active-edge-1} \PP(v_i\sim v_j)= \bigl(1+o(1)
\bigr){\pmatrix{m
\cr
s}}^{-1} \biggl(\E\pmatrix {X_1
\cr
s}
\biggr)^2.
\end{equation}
In particular, we have $ \PP(v_i\sim v_j)=o(1)$.
\end{lem}
%
We note that $\E\bigl({X_i\atop s}\bigr)$ is the expected number of joints
available to the vertex $v_i$.
(\ref{active-edge-1}) may fail in the case where $s=s_m\to\infty$ as
$m\to\infty$; see Example~\ref{ex2}
below.

\begin{remark}\label{rem1}
In the particular case where $s$ remains fixed as
$n,m\to\infty$,
Lemma~\ref{active-edge} suggests that a sparse graph
$G_s(n,m,P)$ is obtained when the random sets are of order
%
\begin{equation}
\label{sparse-D} |D_i|=O_P \bigl(m^{1/2}n^{-1/(2s)}
\bigr).
\end{equation}
By sparse we mean that the degree of the typical vertex is
stochastically bounded as $n\to\infty$.
\end{remark}

Our next result shows that the degree of the typical vertex of a sparse
random intersection graph
is asymptotically
Poissonian
with (random) intensity parameter
$Z_1\E Z_1$. Here $Z_i=Z_{mi}=\bigl({m\atop s}\bigr)^{-1/2}n^{1/2}{\bigl({X_i\atop s}\bigr)}$ denotes the
properly rescaled
number of joints of the vertex $v_i$, $1\le i\le n$.
Let $d(v_i)=d_m(v_i)$ denote the degree of the vertex $v_i$ in $G_s(n,m,P)$.
Observe that, by symmetry, the random variables $d(v_1),\ldots, d(v_n)$
have the same probability distribution.

\begin{tm}\label{Theorem00} Let $m,n\to\infty$. Assume that $s=O(1)$ and:
\begin{longlist}[(iii)]
\item[(i)] $Z_{m1}$ converges in distribution to a random variable $Z$;

\item[(ii)] $\E Z<\infty$ and  $\lim_{m\to\infty}\E Z_{m1}=\E Z$.

Denote $\mu=\E Z$. We have, for $k=0,1,2\ldots,$
%
\begin{equation}
\label{T110} \lim_{m\to\infty} \PP \bigl( d_m(v_1)=k
\bigr)=(k!)^{-1}\E \bigl( (Z\mu )^ke^{-Z\mu} \bigr).
\end{equation}

Formula (\ref{T110}) remains valid if the condition $s=O(1)$ is
replaced by the following weaker condition:
there exists $0<a<1$ such that $s\le am$ and $(s!/n)^{1/s}=o(m)$.
In this case, we require, in addition to \textup{(i)} and \textup{(ii)},
that:

\item[(iii)] $\frac{(s!)^{1/s}}{n^{(s+1)/s}}\sum_{k=2}^n
(Z_{m1}Z_{mk} )^{(s+1)/s}=o_P(1)$.
\end{longlist}
\end{tm}

Formula (\ref{T110}) relates the vertex degree distribution
to the distribution of sizes of the random sets.
In particular, the more variable is the sequence of realized values
$X_1,\ldots, X_n$, the more irregular
is the sequence of vertex degrees of a realized instance
of the intersection graph.
For example,
we obtain a power-law distribution $p_k\sim ck^{-\gamma}$ as $k\to\infty$
whenever the limiting distribution of
$\{Z_{m1}\}_m$ has a power law $\PP(Z>t)\sim c't^{1-\gamma}$ as $t\to \infty$. Here
$\gamma>2$,
and $c,c'$ denote some positive constants.
Observe, that the asymptotic degree distribution defined by formula
(\ref{T110}) has a first moment.

In the case where $s=1$,
the asymptotic degree distribution of an active random intersection graph
was shown (in increasing
generality) in~\cite{stark2004,JKS,Deijfen,Bloznelis2008}; see also~\cite{Bloznelis2010a}
and~\cite{Rybarczyk-degree2011}. In particular,
the result of Theorem
\ref{Theorem00},
for $s=1$, can be found in~\cite{Bloznelis2008,Bloznelis2010a} and
\cite{Rybarczyk-degree2011}. For $s> 1$, the result of Theorem \ref
{Theorem00} is new.
It also applies to the case where $s=s_m\to+\infty$ as $m\to\infty$.
We note that limit (\ref{T110}) may fail when $s=s_m$ grows to infinity
``sufficiently fast;''
see Example~\ref{ex2} below. A general sufficient condition for the convergence
to the
Poisson mixture
(\ref{T110}) is given in Theorem~\ref{Theorem1}.

In Theorem~\ref{Theorem1}, we assume that the random sets $D_{1},\ldots,
D_{n}$ defining the intersection graph
are independent but not necessarily identically distributed. We write,
as before,
$X_{k}=|D_{k}|$ and
assume that for every $k=1,\ldots, n$ and each $A\subset W$, we have $\PP
(D_{k}=A)={\bigl({m\atop  |A|}\bigr)}^{-1}\PP(X_{k}=|A|)$.
Let ${\overline P}$ denote the distribution of the random vector
${\overline X}=(X_{1},\ldots, X_{n})$. By ${\overline
G}_s(n,m,\break{\overline P})$
we denote the random intersection graph
on the vertex
set $V$,
where $v_i, v_j\in V$ are adjacent whenever
$|D_{i}\cap D_{j}|\ge s$. By ${\overline d}(v)={\overline d}_m(v)$ we
denote the degree of
$v\in V$ in ${\overline G}_s(n,m,{\overline P})$.
%
For a vector with nonnegative integer coordinates ${\overline
x}=(x_1,\ldots, x_n)$, we write
$u_k={\bigl({x_1\atop s}\bigr)}{\bigl({x_k\atop s}\bigr)}$ and $x_k^+=\max\{0,x_k-s\}$,
and we denote
\begin{eqnarray*}
\lambda({\overline x})&=&{\pmatrix{m
\cr
s}}^{-1}\sum
_{k=2}^nu_k,\qquad \kappa_1({
\overline x})={\pmatrix{m
\cr
s}}^{-2}\sum_{k=2}^nu_k^2,\\
\kappa_2({\overline x})&=& \frac{x^+_{1}}{m-x_{1}}{\pmatrix{m
\cr
s}}^{-1}\sum_{k=2}^nu_{k}x^+_{k}.
\end{eqnarray*}
%
Our next result applies to a sequence of random intersection graphs
${\overline G}_s(n,m,\break{\overline P})$, where $m\to\infty$ and
$n=n_m\to\infty$, and
where $s=s_m$ and the distribution ${\overline P}$ of the random vector
${\overline X}$ depend on $m$.

\begin{tm}\label{Theorem1} Assume that, as $m\to\infty$:
\begin{longlist}[(vi)]
\item[(iv)] $\lambda({\overline X})$ converges in distribution to a random variable
$\Lambda$;

\item[(v)] $\kappa_1({\overline X})=o_P(1)$;

\item[(vi)] $\kappa_2({\overline X})=o_P(1)$.

Then we have for $k=0,1,2,\ldots,$
%
\begin{equation}
\label{T11} \lim_{m\to\infty} \PP \bigl( {\overline d}_m(v_1)=k
\bigr)=(k!)^{-1}\E \bigl( \Lambda^ke^{-\Lambda} \bigr).\vadjust{\goodbreak}
\end{equation}
\end{longlist}
\end{tm}

We briefly\vspace*{1pt} remark that condition (v)
is
a kind of ``asymptotic negligibility condition'' imposed on the sequence
of random variables ${\bigl({X_1\atop  s}\bigr)}{\bigl({X_k\atop  s}\bigr)}$, $2\le k\le
n$. Condition (vi) ensures that (the distributional limit of) $\lambda
({\overline X})$ alone determines the asymptotic degree distribution;
see also Lemma~\ref{sX1} below.

In Examples~\ref{ex1} and~\ref{ex2} below, we consider \textit{uniform} random
intersection graphs.
In Example~\ref{ex1}, we formulate, in terms of $n=n_m$, $x=x_m$, and $s=s_m$,
conditions that are sufficient for the convergence (\ref{T11}).
Example~\ref{ex2} shows that conditions (iv) and (v) alone do not suffice to establish
the convergence (\ref{T11}). Here we are in a situation
where the edge probability formula (\ref{active-edge-1}) fails.

\begin{example}\label{ex1} Let $\{x_m\}$ and $\{s_m\}$ be integer sequences such
that $1\le x_m-s_m=o ((m-x_m)^{1/2} )$ and
${\bigl({x_m\atop s_m}\bigr)}^2{\bigl({m\atop s_m}\bigr)}^{-1}=o(1)$ as $m\to\infty$. Let
$n_m$ be an integer sequence such that the limit
$\lim_mn_m{\bigl({x_m\atop s_m}\bigr)}^2{\bigl({m\atop  s_m}\bigr)}^{-1}$ exists and is finite. We
denote this limit by $\lambda$.
The sequence of random intersection graphs
$\{G_{s_m}(n_m,   m,  \delta_{x_m})\}$ satisfies the conditions of
Theorem~\ref{Theorem1} with $\Lambda\equiv\lambda$. Therefore, the
degree of the typical vertex has limiting Poisson distribution with
mean $\lambda$.
\end{example}

\begin{example}\label{ex2} Given $m$, let $s=0.5 m$, $x=(\varepsilon+0.5)m$.
Here $\varepsilon\in(0,1)$
is a small absolute constant. We choose $\varepsilon$ small enough
so that $p^*(m):={\bigl({x\atop s}\bigr)}^2{\bigl({m\atop s}\bigr)}^{-1}=o(1)$ as $m\to
\infty$. Next, we choose an
integer sequence
$n_m\uparrow+\infty$
such that $n_m p^*(m)\to1$ as $m\to\infty$ and consider the sequence
of random
intersection graphs $\{G_s(n_m,m, \delta_x)\}_m$.
It follows from the relation
$n_m p^*(m)\to1$ that conditions (iv)--(v) hold with $\Lambda\equiv1$.
Moreover, we show that
${\overline d}_{m}(v_1)=o_P(1)$ as $m\to\infty$.
The latter relation contradicts to (\ref{T11}). Note that for this
sequence of random
intersection graphs, condition
(vi) fails. Indeed, we have
$\kappa_2({\overline X})=\frac{\varepsilon^2m}{0.5-\varepsilon}n_m
p^*(m)\to+\infty$ as
$m\to\infty$. The proofs of the statements of Example~\ref{ex2} are given in
Section~\ref{sec3}.
\end{example}

\subsection{Clustering coefficient}\label{sec2.2}
Typically, the adjacency relations between actors in real networks are
not statistically independent events. Often, chances of a link $v'\sim
v''$ increase as we learn
that actors $v'$ and $v''$ have a
common neighbor,
say,~$v$.
As a theoretical measure of
such a
statistical dependence, one can use the conditional probability
(see, e.g.,~\cite{Deijfen})
\[
\alpha=\PP\bigl(v'\sim v''|
v'\sim v, v''\sim v\bigr).
\]
In the literature (see~\cite{Barrat2000,Newman2001,Newman2003,storgatz1998}), the
empirical estimates of the conditional probability $\alpha$,
\[
{\hat\alpha}=n^{-1}\sum_{v\in V}
\frac{N_3(v)}{N_2(v)}\quad \mbox{and}\quad {\hat{\hat{\alpha}}}=\frac{\sum_{v\in V}N_3(v)}{\sum_{v\in V}N_2(v)},
\]
are called the clustering coefficient and the global clustering
coefficient, respectively.
Here $n$ denotes the number of vertices of a graph, $N_3(v)$ is the
number of unlabeled triangles
having vertex $v$, $N_2(v)$ is the number of unlabeled $2$-stars with
the central vertex $v$.
In this paper, the term clustering coefficient is used exclusively for
the conditional probability
$\alpha$.

Note that, in the random graph $G_s(n,m,P)$, the conditional
probability $\alpha$ does
not depend on the choice of $v,v',v''$. It does not depend on $n$
either. We write
$\alpha=\alpha_s(m,P)$ in order to indicate the dependence on $s, m$
and $P$.

We begin our analysis with the uniform random intersection graph
$G_s(n,\break m,\delta_x)$ where,
for large~$m$,
the asymptotics of the clustering coefficient
is simple and
transparent. Since we are interested in sparse random graphs,
we may assume that $\PP(v_i\sim v_j)=o(1)$ as $m\to\infty$.
Observe that
in the
case of
bounded $s$ [i.e., $s=O(1)$ as $m\to\infty$], this assumption
implies that $x^2=o(m)$; see (\ref{sparse-D}).
In the following lemma, we consider a sequence of uniform random
intersection graphs, where $x=x_m$ and $s=s_m$ depend on $m$,
and $(x-s)^2=o(m-x)$. Note that, for bounded $s$,
the later condition is equivalent to $x^2=o(m)$.

\begin{lem}\label{clustering-c} Let $m\to+\infty$. Assume that $s<x$
and $(x-s)^2=o(m-x)$. Then
%
\begin{equation}
\label{clustering-c-1} \alpha_s(m,\delta_x)={\pmatrix{x
\cr
s}}^{-1}+o(1).
\end{equation}
\end{lem}

The assumption $s<x$ of the lemma excludes the trivial case $x=s$,
where we have $\alpha_s(m,\delta_x)\equiv1$.
Lemma~\ref{clustering-c} brings some insight into the general model
$G_s(n,m,P)$. Namely,
in the case of a sparse random graph,
it suggests
that
the clustering
coefficient is nonvanishing whenever
the sizes of random sets are stochastically bounded as $m\to\infty$.
In fact, we need to impose an even stronger condition, which requires,
in particular, that
the moment of order $2s$ of $P=P_m$
is bounded as
$m\to\infty$.
%
For $k=1,2$, we denote
$a_k=\int{\bigl({x\atop s}\bigr)}^kP(dx)=\E{\bigl({X_1\atop s}\bigr)}^k$.

\begin{lem}\label{Lemma-clustering-1} Assume that for $k=1,2$, we have,
as $m\to\infty$,
%
\begin{equation}
\label{xyze} a_k>0 {\mbox{ and }} \forall \varepsilon\in(0,1)\qquad
\frac{1}{a_k}\int_{x>\varepsilon\sqrt{ m}} {\pmatrix{x
\cr
s}}^k
P(dx)=o(1).
\end{equation}
Then
%
\begin{equation}
\label{clustering-c-2} \alpha_{s}(m,P)=a_1/a_2+o(1).
\end{equation}
\end{lem}

Note that invoking in (\ref{clustering-c-2}) the simple inequality
$a_1^2\le a_2$,
we obtain $\alpha_s(m,P)\le a_2^{-1/2}+o(1)$.
Hence, $a_2\to\infty$ implies $\alpha_s(m,P)\to0$.

Up to our best knowledge, the first result showing that a power-law (active)
random intersection graph with\vadjust{\goodbreak} $m\approx\beta n$, $\beta>0$ admits a
nonvanishing
clustering coefficient is due to Deijfen and Kets \cite
{Deijfen}. They established
a first-order asymptotics of $\alpha_1(m,P)$ as $m\to\infty$
in the particular
case
where $P$ is a mixture of binomial distributions. Yagan and Makowski
\cite{Yagan2009} evaluated the
clustering coefficient $\alpha_1(m,\delta_x)$ of a uniform (active)
random intersection graph and
proved that it is
always positive. Nonvanishing clustering coefficients of a random
intersection digraph
were studied in~\cite{Bloznelis2010}.
The effect of a positive clustering coefficient on the size of the
largest component of $G_1(n,m,P)$ was studied in~\cite{Bloznelis2010a}.
The effect on an epidemic spread was considered in~\cite{Britton2008}.


\subsection{Clustering coefficient and degree}\label{sec2.3}
Here we consider a sequence of sparse random intersection graphs $\{
G_s(n,m,P)\}_m$ with nonvanishing clustering coefficient and
nondegenerate asymptotic vertex degree distribution.

In our next theorem and its corollary, we express an approximate
formula for the clustering coefficient (\ref{clustering-c-2})
in terms of moments
of
the asymptotic vertex degree distribution.
Recall the notation $Z_{1}={\bigl({X_{1}\atop s}\bigr)}{ \bigl({m\atop s}\bigr)}^{-1/2}n^{1/2}$. Here $X_{1}$ is a random
variable with the distribution $P$.

\begin{tm}\label{Theorem-clustering-1+} Let $\beta>0$. Let $m\to\infty$.
Assume that $s=O(1)$ and
%
\begin{equation}
\label{beta+++} {\pmatrix{m
\cr
s}}n^{-1}\to\beta.
\end{equation}
Suppose that conditions \textup{(i)} and \textup{(ii)} of Theorem~\ref{Theorem00} hold.
In the case where:

\textup{(ii$'$)} $0<\E Z^2<\infty$  and  $\lim_{m\to\infty}\E Z_{1}^2=\E Z^2$,

we have
%
\begin{equation}
\label{clustering-c-2+} \alpha_{s}(m,P)=\frac{1}{\sqrt{\beta}}
\frac{\E Z}{\E Z^2}+o(1).
\end{equation}

In the case where $\E Z^2=\infty$, we have $\alpha_s(m,P)=o(1)$.
\end{tm}

\begin{remark}\label{rem2} Observe that
under conditions
(i) and (ii$'$), the asymptotic degree
distribution exists and is defined by (\ref{T110}). Indeed, (i) and
(ii$'$) imply (ii), and,
therefore,
one can apply
Theorem~\ref{Theorem00}.
Furthermore,
$\E Z^2<\infty$
implies that
the asymptotic degree distribution
has a second moment.
Assuming, in addition, that~(\ref{beta+++}) holds, we conclude from
Theorem~\ref{Theorem-clustering-1+}
that the clustering coefficient does not vanish
in the case where the (asymptotic)
degree distribution has finite second moment.
Moreover, if (\ref{beta+++}) fails and we have $n=o\bigl(\bigl({m\atop s}\bigr)\bigr)$, then
the clustering
coefficient vanishes [i.e., $\alpha=o(1)$ as $n,m\to\infty$]
in the case where (i) and (ii$'$) hold, as well as in the case where (i),
(ii) hold and $\E Z^2=\infty$.
The proof of this statement is given in Section~\ref{sec3}.
\end{remark}

Next we express (\ref{clustering-c-2+}) in terms of moments of the
asymptotic degree distribution.
Let $d_*$ be a random variable with the asymptotic degree distribution
defined by (\ref{T110}), that is,
we have
%
\begin{equation}
\label{d-distr} \qquad\PP(d_*=k)=p_k,\qquad p_k=(k!)^{-1}\E
\bigl( (Z\mu)^ke^{-Z\mu} \bigr),\qquad  k=0,1,2,\ldots.\vadjust{\goodbreak}
\end{equation}
Then $\E d_*= (\E Z)^2$ and $\E d_*^2=(\E Z)^2\E Z^2+(\E Z)^2$.
Invoking these formulas
in~(\ref{clustering-c-2+}), we obtain the following:

\begin{col}\label{Corollary-clustering-1+} Let $m,n\to\infty$. Assume
that $s=O(1)$ and conditions
\textup{(i)}, \textup{(ii$'$)} and (\ref{beta+++}) hold. Then we have
%
\begin{equation}
\label{rugs-29} \alpha_{s}(m,P)=\frac{1}{\sqrt{\beta}}\frac{(\E d_*)^{3/2}}{\E d_*^2-\E
d_*}+o(1).
\end{equation}
\end{col}

Foudalis et al.~\cite{Foudalis2011} made an interesting observation,
based on an empirical study of a power-law network data,
that ``clustering correlates negatively with degree.'' More precisely,
their numerical data
suggest that the
conditional probability
%
\begin{equation}
\alpha^{[k]}=\PP\bigl(v'\sim v''
| v\sim v', v\sim v'', d(v)=k\bigr)
\end{equation}
is a decreasing convex function of (the realized value $k$ of) the
degree $d(v)$ of vertex $v$.
Theorem~\ref{Theorem-clustering-degree} establishes a first-order asymptotics
of the conditional probability
$\alpha^{[k]}=\alpha^{[k]}_s(m,P)$ as $m,n\to\infty$ for $G_s(n,m,P)$ model.
The result of Theorem~\ref{Theorem-clustering-degree} provides a rigorous argument explaining
empirical findings of~\cite{Foudalis2011} mentioned above; see also
Example~\ref{ex3} below.

\begin{tm}\label{Theorem-clustering-degree}
Let $\beta>0$. Let $m\to\infty$. Assume that $s=O(1)$ and that
conditions \textup{(i)}, \textup{(ii$'$)} and
(\ref{beta+++}) hold.
Then we have, for every $k=2,3,\ldots,$
%
\begin{equation}
\label{clustering-degree} \alpha^{[k]}_{s}(m,P)=
\frac{1}{k}\frac{\E Z}{\sqrt{\beta}}\frac
{p_{k-1}}{p_{k}}+o(1).
\end{equation}
\end{tm}

Note that, for a power-law asymptotic degree distribution $p_k\sim
ck^{-\gamma}$,
(\ref{clustering-degree})~implies
%
\begin{equation}
\label{k-clustering} \alpha^{[k]}_{s}(m,P)\sim
c'k^{-1}\qquad {\mbox{as }} k\to\infty
\end{equation}
for any $\gamma>3$. Here $c'=\beta^{-1/2}\E Z$.
Hence,
$\alpha^{[\cdot]}$ ``correlates negatively
with degree.''

Next we illustrate (\ref{clustering-degree}) by two examples,
where we assume, for simplicity, that $s\equiv1$.

\begin{example}\label{ex3}
Fix $\gamma>3$ and $\theta>0$, and consider a sequence
of random intersection graphs satisfying the conditions of Theorem \ref
{Theorem-clustering-1+}
with $Z$ having the power-law distribution
$\PP(Z>t)=(\theta t)^{1-\gamma}$ for $t\ge\theta^{-1}$.
In this case, we have $p_k\sim(\gamma-1)\theta^{1-\gamma} k^{-\gamma
}$ as $k\to\infty$. Now
(\ref{clustering-degree}) implies~(\ref{k-clustering}).
\end{example}

\begin{example}\label{ex4}
Fix $\mu>0$ and consider
a sequence
$\{G_1(n,m,\delta_{x})\}$
such that $x(n/m)^{1/2}\to\mu$ as $m\to\infty$. The sequence satisfies
the conditions of
Theorem~\ref{Theorem00}
with $Z$ having the degenerate distribution $\PP(Z=\mu)=1$.
Hence, by Theorem~\ref{Theorem00}, the asymptotic degree distribution is
Poisson with mean $\mu^2$.
In this case, we have $p_{k-1}/p_{k}=k\mu^{-2}$, and (\ref
{clustering-degree}) implies
$\alpha^{[k]}\approx\mu^{-1}$. Hence, $\alpha^{[\cdot]}$ does not correlate
with degree.
\end{example}


The rationale behind these examples is as follows.
In a sparse active graph, an edge
is typically realized by a single joint shared by the two vertices
linked by this edge.
Similarly, a
triangle
is realized by a single joint shared by the three vertices of the triangle.
Therefore, two neighbors, say, $v_i$ and $v_j$ of a given vertex $v_t$
establish a link
whenever the joints responsible for the edges $v_i\sim v_t$ and
$v_j\sim v_t$ match.
In particular, given $X_t=|D_t|$, the probability that two neighbors of
$v_t$ are adjacent is approximately ${\bigl({X_t\atop s}\bigr)}^{-1}$. Here
$\bigl({X_t\atop s}\bigr)$ is the number of joints available to~$v_t$.
%
Next, we remark that given $X_t$, the number of neighbors of $v_t$ has
binomial distribution
$\Bin(n-1, p_{(t)})$
with
$p_{(t)}\approx{\bigl({m\atop s}\bigr)}^{-1}\bigl({X_t\atop s}\bigr)\E\bigl({X_1\atop s}\bigr)$;
see Lemma~\ref{sX1} below.
For large values of $X_t$, this number concentrates around its mean
because of the concentration property of Binomial distribution.
Hence, for large values of $X_t$, the degree $d(v_t)$ scales as
$(n-1)p_{(t)}$ $=\bigl({X_t\atop s}\bigr)\tau$,
where
$\tau\approx\beta^{-1/2}\E Z$ as $n,m\to+\infty$. In particular, given
a vertex $v_t$ with a large degree,
say,
$d(v_t)=k$,
it is reasonable to expect that $(n-1)p_{(t)}\approx k$, that is,
${\bigl({X_t\atop s}\bigr)}^{-1}\approx k^{-1}\tau$.
We summarize the argument as follows: The probability that two
neighbors of a vertex of degree $k$ are adjacent scales
as ${\bigl({X_t\atop s}\bigr)}^{-1}\approx k^{-1}\tau$, for $k\to+\infty$.
This explains formula (\ref{clustering-degree}) in the case where
the sequence
$X_1,X_2,\ldots, X_n$ exhibits a high variability, for example, it is a
sample from a heavy-tailed distribution.
The argument fails in the case of uniform random intersection graphs,
since here all
attribute sets are of the same size ($X_i=\mathrm{const}$) and, hence, the
realized values of $d(v_t)$ and $X_t$
do not correlate.

\begin{remark}\label{rem3}
Using the identity $\E d_*=(\E Z)^2$, we can express (\ref
{clustering-degree}) solely in
terms of the (asymptotic) degree distribution
%
\begin{equation}
\label{clustering-degree++} \alpha^{[k]}_{s}(m,P)=
\frac{1}{k}\frac{\sqrt{\E d_*}}{\sqrt{\beta
}}\frac{\PP(d_*=k-1)}{\PP(d_*=k)}+o(1).
\end{equation}
In the case where $m$ is large and the distribution (\ref{d-distr}) is
close to
that of the
degree sequence of the observed graph, we can replace the moments $\E
d_*, \E d_*^2$ and
probabilities $\PP(d_*=k)$, $\PP(d_*=k-1)$
in (\ref{rugs-29})
and (\ref{clustering-degree++})
by their estimates based on the observed degree sequence. In this way,
we obtain estimates
of $\alpha$ and
$\alpha^{[k]}$ based on
the degree sequence and involving the parameter $\beta$.
\end{remark}

Finally, we note that the second moment $\E Z^2<\infty$ required by
Theorem~\ref{Theorem-clustering-degree}
does not show up in (\ref{clustering-degree}). This observation suggests
indirectly that the second
moment condition could
perhaps be replaced by the weaker first moment condition.

\subsection{Concluding remarks}\label{sec2.4} We remark that the
asymptotic degree distribution of $G_s(n,m,P)$ as $n,m\to\infty$ is
determined by the limiting distribution of the properly scaled number
of joints of a typical vertex, that is, by the limiting distribution of
$Z_1={\bigl({m\atop s}\bigr)}^{-1}n^{1/2}\bigl({X_1\atop s}\bigr)$,
which we denote~$P_Z$. Furthermore, the first-order asymptotics as
$n,m\to+\infty$ of the clustering coefficients $\alpha$ and
$\alpha^{[k]}$ is determined by $P_Z$ and $\beta=\lim\bigl({m\atop
s}\bigr)n^{-1}$.
One may observe that the parameter $s$
does not show up in an explicit way neither in the description of the asymptotic
degree distribution nor in the asymptotic
formulas for $\alpha$ and $\alpha^{[k]}$.
We expect that the role of the parameter~$s$
becomes more important in the case of denser graphs, for example, in
the case where
$\bigl({m\atop s}\bigr)n^{-1}$ converges to zero sufficiently fast.

\section{Passive intersection graph}\label{sec3}

In this section, we consider graphs on the vertex set $W=\{w_1,\ldots,
w_m\}$. Let $s\ge1$ be an integer.
Let $D_1,\ldots, D_n$ be independent random subsets of $W$ having the same
probability distribution~(\ref{W-D-A}). We say that vertices $w,w'\in
W$ are
linked by $D_j$ if $w,w'\in D_j$. For example, every $w'\in D_j\setminus
\{w\}$ is linked to $w$ by $D_j$.
The links created by $D_1,\ldots, D_n$ define a multigraph
on the vertex set $W$.
In the \textit{passive} random
intersection graph, two vertices $w, w'\in W$ are declared adjacent
whenever there are at
least $s$ links between $w$ and $w'$; that is, the pair $\{w,w'\}$
is contained in at least $s$ subsets of the collection $\{D_1,\ldots,
D_n\}$;
see~\cite{godehardt2003}. We denote
the passive random intersection graph $G^*_s(n,m,P)$.
Here $P$ is the common probability distribution of
the random variables $X_1=|D_1|,\ldots, X_n=|D_n|$. We shall consider
only the case where $s=1$.

Before presenting our results, we introduce some notation.
With a sequence of probabilities $Q=\{q_0, q_1,\ldots\}$ such that $\sum_jq_j=1$ and
$\mu_Q=\sum_jjq_j<\infty$,
we associate another sequence
of probabilities ${\tilde Q}=\{{\tilde q}_0, {\tilde q}_1,\ldots\}$
obtained as follows. In the case where
$\mu_Q>0$, we define ${\tilde q}_j=(j+1)q_{j+1}/\mu_Q$, $j=0,1,2,\ldots.$ For $\mu_Q=0$, we put
${\tilde q}_0=1$ and ${\tilde q}_j=0$, $j=1,2,\ldots.$
In particular, we denote by ${\tilde P}={\tilde P}_{X_1}$
the probability distribution on $\{0,1,2,\ldots\}$ putting mass $(j+1)\PP
(X_1=j+1)/\E X_1$ on an
integer $j=0,1,2,\ldots.$
By $P_\xi$ we denote the probability
distribution of a random variable $\xi$.

\subsection{Degree distribution}\label{sec3.1}
Let $d=d(w_1)$ denote the
degree of $w_1$ in $G^*_1(n,\break m, P)$, and let $L=L(w_1)$ denote the number
of links incident to
$w_1$. Our next theorem establishes the asymptotic degree distribution
of a sequence of
sparse passive random intersection graphs
in the case where $m$ and $n$ are of
the same order. By ``sparse'' we mean that the degree $d$ remains
stochastically bounded, that is, $d=O_P(1)$
as $m,n\to\infty$.

\begin{tm}\label{PTh-0} Let $\beta>0$. Let $m,n\to\infty$. Assume that
$mn^{-1}\to\beta$ and:
\begin{longlist}[(viii)]
\item[(vii)] $X_{1}$ converges in distribution to a random variable
$Z$.

\item[(viii)] $\E Z<\infty$ and  $\lim_{m\to\infty}\E X_{1}=\E Z$.\vadjust{\goodbreak}

Then $L$ converges in distribution to the compound Poisson random variable
$d_*=\sum_{j=1}^\Lambda{\tilde Z}_j$.
Here ${\tilde Z}_1,{\tilde Z}_2,\ldots$ are independent random variables with
common probability distribution ${\tilde P}_Z$, the random variable
$\Lambda$
is independent of the sequence
${\tilde Z}_1,{\tilde Z}_2,\ldots$ and has Poisson distribution
with mean $\E\Lambda= \beta^{-1}\E Z$.

If, in addition,

\item[(ix)] $\E Z^{4/3}<\infty$ and  $\lim_{m\to\infty}\E X^{4/3}_{1}=\E Z^{4/3}$,
\end{longlist}
then $d$ converges in distribution to
$d_*$.
\end{tm}

Theorem~\ref{PTh-0} shows how the (asymptotic) vertex degree
distribution depends on the
distribution of sizes of the random sets.
In particular,
we obtain a power-law degree distribution
whenever the distribution
of the sizes
has a power-law.
Moreover, since the distribution of ${\tilde Z}_1$ has a much heavier
tail than that
of $Z$,
we can even obtain (asymptotic) degree distribution with infinite first moment.
We refer to~\cite{Foss} for a survey of results on local and tail probabilities
of sums of heavy-tailed random variables, including, in particular,
random sums and compound Poisson random variables.

\begin{remark}\label{rem4} Note that (vii) and (ix) imply (viii). Hence, condition
(ix) is more restrictive than (viii).
The fact that the asymptotic distribution of $d$ only refers to the
first moment of $Z$ suggests indirectly that the $4/3$ moment condition
(ix) could perhaps be
waived.
\end{remark}

In the next lemma, we collect several facts about the degree distribution
in the case where either $n=o(m)$ or $m=o(n)$.

\begin{lem}\label{PLemma-0}
Let $m,n\to\infty$. In the case where $n=o(m)$, we have
%
\begin{equation}
\label{m>>n} \forall \varepsilon>0 \qquad\PP \bigl(d\ge\varepsilon m n^{-1}-1
\bigr)\ge\PP(d\ge1)-\varepsilon.
\end{equation}

Now assume that $m=o(n)$. Denote $n_*=n\PP(X_1\ge2)$. We distinguish
three cases:
\textup{(a)} $n_*=o(m)$;
\textup{(b)} $m=o(n_*)$;
\textup{(c)} $mn_*^{-1}=\beta(1+o(1))$ for some $\beta>0$.
In the case \textup{(a)}, we have
%
\begin{equation}
\label{n*} \forall \varepsilon>0 \qquad \PP \bigl(d\ge\varepsilon m \bigl(\max\{1,n_*\}
\bigr)^{-1}-1 \bigr)\ge\PP(d\ge 1)-\varepsilon.
\end{equation}
In the case \textup{(b)}, we have
%
\begin{equation}
\label{case-b} \forall C >0 \qquad \PP(d>C)=1-o(1).
\end{equation}
In the case \textup{(c)}, the conclusion of
Theorem~\ref{PTh-0} holds if the conditions of the theorem are
satisfied with the random variable
$X_1$ replaced by the random variable $X_{1*}$ having the distribution
$\PP(X_{1*}=i)=\PP(X_1=i|X_1\ge2)$, $i\ge2$.
\end{lem}

\begin{remark}\label{rem5} Observe that any of inequalities
(\ref{m>>n}), (\ref{n*}) and (\ref{case-b})
rules out the
option of a nondegenerate and\vadjust{\goodbreak} stochastically bounded $d$ as $m,n\to
\infty$.
%
Therefore, a nontrivial asymptotic degree distribution is possible only
in the cases where
$c_1<m^{-1}n<c_2$ or
$c_1<m^{-1}n\PP(X_1\ge2)<c_2$ as $n,m\to\infty$ for some $c_1,c_2>0$.
Furthermore, in the second case,
the number $N=\sum_{i=1}^n{\mathbb I}_{\{X_i\ge2\}}$ of random sets
that define the
edges of $G^*_1(n,m,P)$ concentrates around the expected value $n_*=\E
N=n\PP(X_1\ge2)$, and, for this reason,
the asymptotic vertex degree distribution of $G^*_1(n,m,P)$
is the same as that of $G^*_1(\lfloor n_*\rfloor,m,P_*)$. Here
$\lfloor n_*\rfloor$ denotes the largest integer not exceeding
$n_*$,
and $P_*$ denotes the distribution of $X_{1*}$.
Let us note that the asymptotic degree distribution of $G^*_1(\lfloor
n_*\rfloor,m,P_*)$ can be obtained
from Theorem~\ref{PTh-0}.
\end{remark}

The asymptotic degree distribution of the passive random intersection
graph $G^*_1(n,m,P)$ has been studied by Jaworski and Stark
\cite{JaworskiStark2008}.
They showed a necessary and sufficient condition for the convergence of
the degree distribution
and determined conditions for the convergence to a Poisson limit.
They also asked what other possible limiting distributions are.
Theorem~\ref{PTh-0} and Lemma
\ref{PLemma-0} answer this question in the particular case of sparse graphs.
Let us mention that the approach used in the proof of Theorem \ref
{PTh-0} is different from that of~\cite{JaworskiStark2008}.

\subsection{Clustering coefficient and degree}\label{sec3.2}
The
clustering coefficient
\[
\alpha^*=\alpha^*(n,m,P)=\PP ( w_2\sim w_3 |
w_1\sim w_2, w_1\sim w_3 )
\]
of a passive random intersection graph $G^*_1(n,m,P)$ has been studied
in the recent paper
by Godehardt et al.~\cite{GJR2010}.
They showed, in particular, that
%
\begin{equation}
\label{alpha-*} \qquad\alpha^*(n,m,P)=\frac{\beta_*^2m^{-1} ( \E(X_1)_2  )^3+ \E
(X_1)_3} {\beta_* ( \E(X_1)_2  )^2+\E(X_1)_3}+o(1),\qquad \beta_*:=nm^{-1},
\end{equation}
provided that $\E(X_1)_2>0$ and $\E(X_1)_2 =o(m^2n^{-1})$ as $m,n\to
\infty$; see Theorem 2 and Corollary~1 in~\cite{GJR2010}. Here and
below we denote $(x)_k=x(x-1)\cdots(x-k+1)$.

We are interested in the relation between the clustering coefficient
and the degree.
Our first result expresses (\ref{alpha-*}) in terms of moments of the
(asymptotic) degree distribution.


\begin{tm}\label{PTh-1} Let $\beta>0$. Let $m,n\to\infty$ so that
$mn^{-1}\to\beta$.
Suppose that condition \textup{(vii)} of Theorem~\ref{PTh-0} holds and that, in addition,
$0<\E Z^2<\infty$ and  $\lim_{m\to\infty}\E X^2_{1}=\E Z^2$.
Then the degree $d$ converges in distribution to the compound Poisson random
variable $d_*$ defined in Theorem~\ref{PTh-0}.

In the case where $\PP(Z\ge2)>0$, $\E Z^3<\infty$
and  $\lim_{m\to\infty}\E X^3_{1}=\E Z^3$, we have
%
\begin{equation}
\label{alpha-*+} \alpha^*(n,m,P)=\frac{ \E(d_*)_2-(\E d_*)^2}{\E(d_*)_2}+o(1).
\end{equation}
In the case $\E Z^3=\infty$, we have $\alpha^*(n,m,P)=1+o(1)$.
\end{tm}

Observe that $\E Z^3<\infty$ if and only if $\E d_*^2<\infty$.
Hence, in order to obtain the clustering coefficient $\alpha^*<1$ (as
$n,m\to\infty$), we need to require
$\E d_*^2<\infty$.

\begin{remark}\label{rem6} The result of Theorem~\ref{PTh-1} extends to the case
where $m=o(n)$ and
$mn_*^{-1}\to\beta$. Indeed, in this case, one can apply
Theorem~\ref{PTh-1} to the random graph
$G^*_1(\lfloor n_*\rfloor, m ,P_*)$; see Remark~\ref{rem5}.
\end{remark}

Next, we examine the conditional probability
\begin{eqnarray*}
\alpha^{*[k]} = \alpha^{*[k]}(n,m,P) = \PP \bigl(
w_2\sim w_3 | w_1\sim w_2,
w_1\sim w_3, d(w_1)=k \bigr).
\end{eqnarray*}

\begin{tm}\label{PTh-2}
Let $\beta>0$. Let $m,n\to\infty$ so that $mn^{-1}\to\beta$.
Assume that condition \textup{(vii)} of Theorem~\ref{PTh-0} holds.
Assume, in addition, that
the random variable~$Z$ defined by \textup{(vii)}
has the third moment, $\PP(Z\ge2)>0$ and  $\lim_{m\to\infty}\E
X^3_{1}=\E Z^3$.
Denote $d_{2*}=\sum_{j=1}^{\Lambda}({\tilde Z}_j)_2$, where
$\Lambda, {\tilde Z}_1,{\tilde Z}_2,\ldots$ are defined in
Theorem~\ref{PTh-0}. Then for every $k=2,3,\ldots$ satisfying $\PP
(d_*=k)>0$, we have
%
\begin{equation}
\label{alpha-*++} \alpha^{*[k]}(n,m,P)=\frac{1}{k(k-1)} \E (
d_{2*} | d_*=k )+o(1).
\end{equation}
\end{tm}

\begin{example}\label{ex5} Given $\beta>0$ and an integer $x\ge2$, consider the
random graph $G^*_1(n,m,\delta_x)$, where $m=\lfloor\beta n\rfloor$ as
$m,n\to\infty$.
By Theorem~\ref{PTh-0}, the degree $d$ converges in distribution
to the random variable $d_*=(x-1)\Lambda$. Furthermore, we have
$d_{2*}=(x-1)_2\Lambda$.
Here the random variable $\Lambda$ has Poisson distribution with mean
$\beta^{-1}(x-1)$.
It follows from (\ref{alpha-*++}) that for every $k=t(x-1)$,
$t=1,2,\ldots,$ we have
%
\begin{equation}
\label{delta-x} \alpha^{*[k]}(n,m,\delta_x) 
=
\frac{x-2}{k-1}+o(1) \qquad{\mbox{as }} m,n\to\infty.
\end{equation}
Hence, in this case, $\alpha^{*[k]}$ is of order $k^{-1}$ as $k\to\infty
$. Therefore,
$\alpha^{*[\cdot]}$ ``correlates negatively
with degree.''
\end{example}

In order to explain (\ref{delta-x}), we first note that a triangle in a
sparse passive
random intersection graph is typically realized by a single set $D_j$
that covers all
three vertices. Furthermore, with a high probability, any two sets
$D_i, D_j$
that cover a given vertex $w_1$ have no other element in common, that
is, $D_i\cap D_j=w_1$.
Therefore, almost all triangles incident to $w_1$ are realized by
the sets $D_{i_1},\ldots, D_{i_t}$ that cover $w_1$, and the number of such
triangles $N_{\Delta}$ is approximately the sum of numbers
of pairs $\{w', w''\}$ covered by $D_{ij}$, $1\le j\le t$.
Hence, $N_{\Delta}\approx\bigl({|D_{i_1}|-1\atop  2}\bigr)+\cdots+\bigl({|D_{i_t}|-1\atop 2}\bigr)$.
In addition, the degree $d(w_1)\approx|D_{i_1}|-1+\cdots+|D_{i_t}|-1$.
Hence,
given $D_{i_1},\ldots, D_{i_t}$,
the probability that a pair of neighbors of $w_1$
makes a triangle with $w_1$ approximately equals $N_{\Delta}/
\bigl({d(w_1)\atop 2}\bigr)$.\vadjust{\goodbreak}
In the case where the sizes of random sets do not deviate much from
their average value
$\E X_1$, we write $|D_{ij}|\approx\E X_1$ and obtain
$N_{\Delta}/\bigl({d(w_1)\atop 2}\bigr)\approx
(\E X_1-2)/(d(w_1)-1)$. This explains and generalizes~(\ref{delta-x}).

In the case where sizes of random sets are heavy-tailed random variables,
we expect a different pattern for large values of $d(w_1)$. Now the
maximal size
$\max_{1\le j\le t}|D_{ij}|$
is of the same order as the sum $\sum_{1\le j\le t} |D_{ij}|$.
For this reason, the fraction $N_{\Delta}/\bigl({d(w_1)\atop 2}\bigr)$ stays
bounded away from zero.
Moreover, one may expect that, in this case, $\alpha^{*[k]}\to1$ as
$k\to\infty$.


\section{Proofs}\label{sec4}

We begin with general lemmas that are used in the proofs below. Then we
prove results for the active
graph.
Afterwards, we prove results for the passive graph. We note that the
notation introduced in the proof of a
particular lemma or theorem is only valid for that proof.

\subsection{General lemmas}\label{sec4.1}
The following inequality is referred to as LeCam's lemma; see, for
example,~\cite{Steele}.
%
\begin{lem}\label{LeCam1} Let $S={\mathbb I}_1+ {\mathbb I}_2+\cdots+
{\mathbb I}_n$ be the sum of independent random indicators
with probabilities $\PP({\mathbb I}_i=1)=p_i$. Let $\Lambda$ be Poisson
random variable with mean $p_1+\cdots+p_n$. The total variation
distance between the distributions $P_S$ and $P_{\Lambda}$ of $S$ and
$\Lambda$ satisfies the inequality
%
\begin{equation}
\label{LeCam} d_{\mathrm{TV}}(P_{S},P_{\Lambda})\le2\sum
_{i}p_i^2.
\end{equation}
\end{lem}

\begin{lem}\label{sX1} Given integers $1\le s\le d_1\le d_2\le m$, let
$D_1,D_2$ be independent random subsets of the set $W=\{1,\ldots, m\}$
such that
$D_1$ (resp., $D_2$) is uniformly distributed in the class of subsets of
$W$ of size $d_1$ (resp., $d_2$).
The probabilities
${\mathring{p}}:=\PP(|D_1\cap D_2|=s)$
and
${\tilde p}:=\PP(|D_1\cap D_2|\ge s)$ satisfy
%
\begin{equation}
\label{sp} \biggl(1-\frac{(d_1-s)(d_2-s)}{m+1-d_1} \biggr)p^*_{d_1,d_2,s} \le {
\mathring{p}} \le {\tilde p} \le p^*_{d_1,d_2,s},
\end{equation}
where we denote $p^*_{d_1,d_2,s}={\bigl({d_1\atop  s}\bigr)}{\bigl({d_2\atop  s}\bigr)}
{\bigl({m\atop  s}\bigr)}^{-1}$.
\end{lem}
\begin{pf}
It suffices to establish
inequalities (\ref{sp})
for conditional proba\-bilities given~$D_2$.
In order to prove the left inequality, we write
${\mathring{p}}=\break{\bigl({d_2\atop  s}\bigr)}{\bigl({m-d_2\atop  d_1-s}\bigr)}{\bigl({m\atop d_1}\bigr)}^{-1}=y  p^*_{d_1,d_2,s}$,
where
%
\begin{eqnarray*}
y= \prod_{i=0}^{d_1-s-1} \biggl( 1-
\frac{d_2-s}{m-s-i} \biggr) \ge 1- \sum_{i=0}^{d_1-s-1}
\frac{d_2-s}{m-s-i} \ge 1-\frac{(d_2-s)(d_1-s)}{m+1-d_1}.
\end{eqnarray*}
Let us show the right inequality of (\ref{sp}). Since, for every
$D\subset D_2$ of size $|D|=s$, the number of
subsets of $W$ of size $d_1$ that contain $D$ is at most $\bigl({m-s\atop  d_1-s}\bigr)$, we conclude
that ${\tilde p}\le{\bigl({d_2\atop  s}\bigr)}{\bigl({m-s\atop  d_1-s}\bigr)}{\bigl({m\atop  d_1}\bigr)}^{-1}$.
Note that the quantity in the right-hand side of the latter inequality
equals
$p^*_{d_1,d_2,s}$.
\end{pf}

\begin{lem}\label{auxiliarylema} Let $\{X_{n1},X_{n2},\ldots, X_{nn}\}_{n\ge1}$ be a
collection of nonnegative random variables
such that,
for each $n$, the random variables $X_{n1},\ldots, X_{nn}$ are
independent and identically distributed.
Let $\alpha\ge1$. Let $S_{\alpha}=n^{-\alpha}(X_{n1}^{\alpha}+\cdots
+X_{nn}^{\alpha})$.
Assume that $\E X_{n1}<\infty$ for each $n$ and
that the sequence $\{X_{n1}\}_n$ converges in distribution to a random
variable $Z$. Assume, in addition, that $\E Z<\infty$ and
$\lim_n\E X_{n1}=\E Z$.
Then
%
\begin{eqnarray}
\label{FirstFact}  \sup_n\E X_{n1}{\mathbb
I}_{\{X_{n1}>x\}}&\to&0 \qquad{\mbox{as }} x\to+\infty,
\\
\label{SecondFact}  S_1-\E Z&=&o_P(1)\qquad {\mbox{as }} n
\to+\infty,
\\
\label{ThirdFact} \forall \alpha>1 \qquad S_\alpha&=&o_P(1) \qquad{\mbox{as }} n\to+\infty.
\end{eqnarray}
\end{lem}
\begin{pf}
Relations (\ref{FirstFact}) and (\ref{SecondFact}) are shown in the
proof of
Remark 1 in~\cite{Bloznelis2008} and Corollary 1 in~\cite{Bloznelis2010}.
Let us prove (\ref{ThirdFact}). We shall show that
$\mathbf{P}(S_\alpha > \varepsilon) \leq \varepsilon \mathbf{E}Z + o(1)$
as $n\to\infty$
for each $\varepsilon\in(0,1)$.
Given $\varepsilon$, introduce the event
$\mathcal{ B}=\{\max_{1\le i\le n}X_{ni}<\varepsilon^{2/(\alpha-1)}n\}$
and note that $S_{\alpha}\le\varepsilon^2S_1$
when the event $\mathcal{ B}$ holds.
Hence, we have
\begin{eqnarray*}
\PP(S_{\alpha}>\varepsilon) \le \PP \bigl(\{S_{\alpha}>\varepsilon\}
\cap\mathcal{ B} \bigr) + \PP({\overline \mathcal{ B}}) \le\PP
\bigl(S_1>\varepsilon^{-1}\bigr)+\PP({\overline\mathcal{
B}}).
\end{eqnarray*}
Here the complement
event ${\overline\mathcal{ B}}$ has the probability
\begin{eqnarray*}
\PP({\overline\mathcal{ B}}) \le n\PP\bigl(X_{n1}\ge
\varepsilon^{2/(\alpha-1)}n\bigr) \le \varepsilon^{-2/(\alpha-1)}\E
X_{n1}{\mathbb I}_{\{X_{n1}\ge\varepsilon
^{2/(\alpha-1})n\}} =o(1). 
%
\end{eqnarray*}
In the last step, we applied (\ref{FirstFact}). Finally, the bound
$\PP(S_1>\varepsilon^{-1})\le\varepsilon\E S_1=\varepsilon(\E Z+o(1))$ completes the proof.
\end{pf}
%
%

\subsection{Active graph}\label{sec4.2}
We first note that Lemma~\ref{active-edge} is an immediate consequence
of Lemma~\ref{sX1}. Next, we prove results on the degree distribution:
Theorems~\ref{Theorem1} and~\ref{Theorem00} and Example~\ref{ex2}. Afterwards, we prove the statements
related to clustering coefficient.


\begin{pf*}{Proof of Theorem~\ref{Theorem1}} Let $d_{\infty}$ be a
random variable with the
distribution
defined by the right-hand side of (\ref{T11}). We write, for short,
$d_m={\overline d}_m(v_1)$.
Let $f_m(t)=\E e^{it d_m}$
and $f(t)=\E e^{itd_{\infty}}$ denote the
Fourier transforms of the probability distributions of
$d_m$
and $d_{\infty}$. In order to prove the theorem,
we show
that $\lim_mf_m(t)=f(t)$ for every real $t$.

Given $0<\delta<0.01$ and an integer $m$, introduce the event
$\mathcal{ A}=\{\kappa_i({\overline X})<\delta,  i=1,2\}$. Note that,
by (v) and (vi), we have
$\PP(\mathcal{ A})=1-o(1)$ as $m\to\infty$.
Therefore, we write
%
\begin{equation}
\label{f-1} f_m(t)=\E \bigl( e^{itd_m}{\mathbb
I}_{\mathcal{ A}} \bigr) +o(1).
\end{equation}
%
On the event $\mathcal{ A}$, we approximate the conditional characteristic function
\begin{eqnarray*}
f_{m}(t;\overline{x}):=\E \bigl(e^{itd_m} |{\overline X}={
\overline x} \bigr) 
\end{eqnarray*}
by the Fourier transform of the
Poisson distribution with mean $\lambda(\overline{x})$,
\[
g_m(t;\overline{x})=\exp\bigl\{\lambda(\overline{x})
\bigl(e^{it}-1\bigr)\bigr\}.
\]
Since the conditional distribution of $d_m$, given the event $\{
{\overline X}={\overline x}\}$,
is that of the sum of independent Bernoulli random variables with
success probabilities
\begin{eqnarray*}
q_k=\PP \bigl(|D_{1}\cap D_{k}|\ge s
|X_{1}=x_1, X_{k}=x_k \bigr),\qquad  2\le k
\le n,
\end{eqnarray*}
we write
\begin{eqnarray*}
f_{m}(t;\overline{x}) = \prod_{2\le k\le n}
\bigl(1+q_k\bigl(e^{it}-1\bigr) \bigr) = \exp \biggl\{ \sum
_{2\le k\le n}\ln \bigl(1+q_k
\bigl(e^{it}-1\bigr) \bigr) \biggr\}.
\end{eqnarray*}
Note that (\ref{sp})
implies
%
\begin{equation}
\label{q-2i} {\pmatrix{m
\cr
s}}^{-1}u_k \biggl(1-
\frac{x_1^+x_k^+}{m-x_1} \biggr) \le q_k \le {\pmatrix{m
\cr
s}}^{-1}u_k,\qquad 2\le k\le n.
\end{equation}
It follows from the right-hand side inequality of (\ref{q-2i}) and the
inequality
${\bigl({m\atop  s}\bigr)}^{-2}u_k^2\le\delta\le0.01$, which holds on the event
$\mathcal{ A}$,
that
for each $k$, we have $ |q_k(e^{it}-1) |<0.5$.
Invoking the inequality $|\ln(1+z)-z|\le|z|^2$ for complex numbers $z$
satisfying $|z|\le0.5$
(see, e.g., Proposition 8.46 of~\cite{Breiman}), we obtain
from (\ref{q-2i}) that
\begin{eqnarray*}
f_{m}(t;\overline{x})=\exp \bigl\{ \lambda(\overline {x})
\bigl(e^{it}-1\bigr)+r(t) \bigr\},
\end{eqnarray*}
where $|r(t)|\le4 \kappa_1(\overline{x})+2\kappa_2(\overline{x})$.
Now, the inequalities
$\kappa_i(\overline{x})<\delta$, $i=1,2$, which hold on the event
$\mathcal{ A}$,
imply that
$ |f_{m}(t;\overline{x})-g_m(t;\overline{x}) |\le7\delta$.
Invoking this inequality in (\ref{f-1}), we obtain
\begin{eqnarray*}
\bigl|f_m(t)-\E\exp \bigl\{\lambda({\overline X}) \bigl(e^{it}-1
\bigr) \bigr\}\bigr |\le7\delta+o(1)\qquad {\mbox{as }} n\to\infty.
\end{eqnarray*}
Finally, the convergence in distribution of $\{\lambda({\overline X})\}
$ [i.e., condition (iv)]
implies the convergence of the corresponding expectations
of bounded continuous functions. Therefore, $\lim_m\E e^{\lambda
({\overline X})(e^{it}-1)}=f(t)$.
We obtain the inequality $\limsup_m|f_m(t)-f(t)|\le7\delta$, which
holds for
arbitrarily small $\delta>0$.
The proof of Theorem~\ref{Theorem1} is complete.
\end{pf*}

\begin{pf*}{Proof of Theorem~\ref{Theorem00}} Throughout the proof,
we assume that
the integers $n=n_m$, $s=s_m$, and the distributions of random variables
$(X_1,\ldots,\break X_n)={\overline X}$ and
$Z_i=\bigl({X_i\atop  s}\bigr)n^{1/2}{\bigl({m\atop s}\bigr)}^{-1/2}$, $1\le i\le n$, all
depend on $m$.

We derive Theorem~\ref{Theorem00} from Theorem~\ref{Theorem1}.
To this aim, we verify conditions~(v) and (vi) and show that condition
(iv) holds with $\Lambda=\mu Z$.
The fact that (i) and (ii) imply (iv) and (v) follows from Lemma \ref
{auxiliarylema}. Here we show that
(i), (ii) and (iii) imply (vi) in the case where $s\le am$ and
$(s!/n)^{1/s}=o(m)$. In addition,
we show that (i) and (ii) imply (iii) in the case where $s=O(1)$. Note that
$s=O(1)$ means that the sequence $\{s_m\}$ is bounded, and this is a
much more restrictive
condition than $s\le am$ and $(s!/n)^{1/s}=o(m)$.

We start with the observation that
the inequality $1-y^{-1}\le\sqrt{1-m^{-1}}$, which
holds for
$1\le y\le m$, implies the inequality
%
\begin{equation}
\label{formule} \frac{{\bigl({x\atop s}\bigr)}}{\sqrt{\bigl({m\atop s}\bigr)}} = \frac{(x)_s}{\sqrt{(m)_s}} \frac{1}{\sqrt{s!}} \ge
\biggl( \frac{x-s+1}{\sqrt{m-s+1}} \biggr)^s \frac{1}{\sqrt{s!}} \qquad\forall x
\in[s,m].
\end{equation}

In order to show (vi), we prove that
%
\begin{equation}
\label{kappa+1+2} \forall\varepsilon\in(0,1)\qquad \limsup\PP \bigl(
\kappa_2({\overline X})>\varepsilon \bigr)\le \varepsilon.
\end{equation}
Here and below, the limits are taken as $m\to\infty$.
Given $\varepsilon$, we find a (sufficiently large) constant $A>0$ such that
$\PP(Z_{1}\ge A)\le A^{-1}\E Z_{1}<\varepsilon$
uniformly in $m$.
Here we applied Markov's inequality
and used (ii). In view of the inequality $\PP(Z_{1}\ge A)<\varepsilon$,
we can write
%
\begin{equation}
\label{kappa-kappa} \PP \bigl( \kappa_2({\overline X})>\varepsilon
\bigr)\le\varepsilon +p_{\varepsilon},
\end{equation}
where
$p_{\varepsilon}:=\PP ( \{\kappa_2({\overline X})>\varepsilon\}\cap
\{ Z_{1}<A\} )$.
Observing that on the event $\{Z_{1}<A\}$ we have, by (\ref{formule}),
\[
X_1-s+1\le\sqrt{m-s+1} \Bigl(A\sqrt{s!/n} \Bigr)^{1/s},
\]
we obtain from the bound $(s!/n)^{1/s}=o(m)$ the inequality
$X_1\le am+o(m)$.
It follows now that
\[
Y_i:=\frac{(X_1-s+1)(X_i-s+1)}{m-X_1+1} \le \bigl((1-a)^{-1}+o(1)
\bigr)\frac{(X_1-s+1)(X_i-s+1)}{m}.
\]
Furthermore, invoking the inequality
%
\begin{equation}
\frac{(X_1-s+1)(X_i-s+1)}{m} \le \biggl( s! \frac{ \bigl({X_1\atop  s}\bigr) \bigl({X_i\atop  s}\bigr) }{ \bigl({m\atop  s}\bigr) } \biggr)^{1/s},
\end{equation}
which follows from (\ref{formule}),
we obtain the inequality
\begin{eqnarray*}
Y_i\le \bigl((1-a)^{-1}+o(1) \bigr) \bigl(\bigl(Z_{1}Z_{i}\bigr)s!n^{-1}
\bigr)^{1/s}.
\end{eqnarray*}
Hence, on the event $\{Z_{1}\le A\}$, we have $\kappa_2({\overline
X})\le ((1-a)^{-1}+o(1) )\varkappa$.
Here we denote $\varkappa=\frac{({s!})^{1/s}}{n^{(s+1)/s}}\sum_{k=2}^nt_k^{1+s}$ and $t_k= (Z_{1}Z_{k})^{1/s}$.
Now we see that
(iii) implies $p_{\varepsilon}=o(1)$, and, therefore, (\ref{kappa+1+2})
follows from (\ref{kappa-kappa}).

In the remaining part of the proof, we show that (i) and (ii), together
with the condition $s=O(1)$,
imply (iii). For $s\ge2$, we write, by
H\"older's inequality,
%
\begin{equation}
\label{Hoelder} \frac{n^{(s+1)/s}}{({s!})^{1/s}}\varkappa=\sum_{k=2}^n
t_k^{s-1} t_k^{2}\le \Biggl(\sum
_{k=2}^nt_k^{s}
\Biggr)^{(s-1)/s} \Biggl(\sum_{k=2}^nt_k^{2s}
\Biggr)^{1/s}
\end{equation}
and observe that (\ref{Hoelder}) implies
\begin{eqnarray*}
(s!)^{-1/s}\varkappa \le \bigl(\lambda({\overline X})
\bigr)^{(s-1)/s} \bigl(\kappa_1({\overline X})
\bigr)^{1/s}.
\end{eqnarray*}
For $s=1$, we have $\varkappa=\kappa_1({\overline X})$.
Now (iv) and (v), which follow from (i) and (ii), imply that $\varkappa
=o_P(1)$ as $m\to\infty$.
\end{pf*}

\begin{pf*}{Proof of Example~\ref{ex2}}
Let $D_1$, $D_2$ be independent random subsets which are uniformly
distributed in the class of subsets of
$[m]=\{1,2,\ldots, m\}$ of size $x=(\varepsilon+0.5)m$. Let $s=0.5 m$. Denote
\begin{eqnarray*}
p^*(m)&=&{\pmatrix{x
\cr
s}}^2{\pmatrix{m
\cr
s}}^{-1},\qquad
p'(m)=\PP\bigl(|D_1\cap D_2|= s\bigr),\\
p''(m)&=&\PP\bigl(|D_1\cap D_2|\ge
s\bigr).
\end{eqnarray*}
We show that for a sufficiently small absolute constant $\varepsilon\in
(0,1)$, we have, as $m\to\infty$,
%
\begin{equation}
\label{pp23} p^*(m)=o(1),\qquad  p''(m)=o \bigl(p^*(m)
\bigr).
\end{equation}
Observe that, by the first relation of (\ref{pp23}), we can construct
an increasing integer sequence
$\{n_m\}$ such that $n_mp^*(m)\to1$ and
$n_m (p^*(m))^2\to0$ as $m\to\infty$. Hence, conditions (iv) and (v) of
Theorem~\ref{Theorem1} are fulfilled with $\Lambda\equiv1$.
In addition, by the second relation of (\ref{pp23}), the average degree
$\E  d_{m}(v_1)=(n_m-1)p''(m)$ satisfies $\E  d_{m}(v_1)=o(1)$.
Hence, $d_{m}(v_1)=o_P(1)$, and this means that (\ref{T11}) fails.

Let us prove (\ref{pp23}). By Stirling's formula,
as $m\to\infty$, we have $\bigl({m\atop s}\bigr)\sim2^m(0.5\pi m)^{-1/2}$ and
\[
{\pmatrix{x
\cr
s}}^2\sim \biggl(\frac{x}{x-s}
\biggr)^{2(x-s)} \biggl(\frac
{x}{s} \biggr)^{2s}
\frac{x}{2\pi s(x-s)} = \frac{0.5+\varepsilon}{\pi\varepsilon m} A_{\varepsilon}^m,
\]
where $A_{\varepsilon}= ( 1+(2\varepsilon)^{-1} )^{2\varepsilon
}(1+2\varepsilon)\to1$ as $\varepsilon\to0$.
In particular, for a small $\varepsilon\in(0,1)$, we have
$A_{\varepsilon}<2$. The latter inequality implies $p^*(m)=o(1)$.

Let us prove the second relation of (\ref{pp23}). We denote $p'_r=\PP
(|D_1\cap D_2|=r )$ and write
%
\begin{equation}
\label{B-epsilon} p''(m)=\sum
_{r=s}^{x}p'_r=p'_s
\Biggl(1+\sum_{k=1}^{x-s}
\bigl(p'_{s+k}/p'_s\bigr)
\Biggr).
\end{equation}
%
It follows from the identities
$p'_r=\bigl({x\atop r}\bigr)\bigl({m-x\atop  x-r}\bigr){\bigl({m\atop  x}\bigr)}^{-1}$ and
$\frac{p'_{s+k}}{p'_s}
=
\frac{(x-s)_k  (x-s)_k}{(s+k)_k  (m-2x+s+k)_k}$
that
%
\begin{eqnarray*}
\frac{p'_{s+k}}{p'_s}\le\frac{(x-s)^{2k}}{s^k(m-2x+s)^k}= B_{\varepsilon}^k,
\end{eqnarray*}
where $B_\varepsilon= 4\varepsilon^2/(1-4\varepsilon)\to0$ as
$\varepsilon\to0$.
In particular, we have $B_{\varepsilon} <0.07$ for $\varepsilon<0.1$.
Invoking these inequalities in (\ref{B-epsilon})\vadjust{\goodbreak} and observing that
$p'_s=p'(m)$, we obtain,
for $\varepsilon<0.1$,
\begin{eqnarray*}
p'(m)<p''(m)< 1.1
p'(m).
\end{eqnarray*}
We complete the proof of (\ref{pp23}) by showing that $p'(m)=o
(p^*(m) )$. We have
\[
\hspace*{-6pt}\frac{p'(m)}{p^*(m)}=\frac{(m-x)_{x-s}}{(m-s)_{x-s}}\le \biggl(\frac
{m-x}{m-s}
\biggr)^{x-s} = (1-2\varepsilon )^{\varepsilon m}=o(1) \qquad{\mbox{as }} m\to
\infty.\hspace*{10pt}
\]
\upqed\end{pf*}

%

Before proving the statements related to clustering coefficient, we
introduce some notation. For a real number $a$, we denote $a_+=\max\{0,
a\}$.
Recall that $V=\{v_1,\ldots, v_n\}$ is the vertex set of $G_s(n,m,P)$,
and $X_i=|D_i|$ denotes the size of the attribute set $D_i$
of $v_i$, $1\le i\le n$.
By ${\tilde\PP}(\cdot)=\PP(\cdot|D_1, X_2, \ldots, X_n)$ we denote the
conditional probability given $D_1$ and $X_2,\ldots, X_n$.
Write, for short, $M=\bigl({m\atop  s}\bigr)$ and denote $Y_i=\bigl({X_i\atop  s}\bigr)$,
$Z_i=Y_in^{1/2}M^{-1/2}$, and $a_k=\E Y_1^k$.
We write $D_{ij}=D_i\cap D_j$ and
introduce the events $\mathcal{E}_{ij }= \{v_i \sim v_j\}$,
\begin{eqnarray*}
\mathcal{ A}_1&=&\{v_1\sim v_2,
v_1\sim v_3, v_2\sim v_3\},\qquad
\mathcal{ A}_2=\bigl\{|D_{12}|=|D_{13}|=|D_{23}|=s
\bigr\},
\\
\mathcal{ A}_3&=&\{D_{12}=D_{13}=D_{23}
\}\cap\mathcal{ A}_2,\qquad \mathcal{ B}=\{v_1\sim
v_2, v_1\sim v_3\},
\\
 \mathcal{ C}&=& \bigl\{|D_{23}\setminus D_1|\ge1,
|D_{12}|=|D_{13}|=s\bigr\}, \qquad \mathcal{ D}_{ij}=
\bigl\{|D_{ij}|\ge s+1\bigr\}.
\end{eqnarray*}
We have
%
\begin{equation}
\label{ALPHA-H1} \alpha=\alpha_s(m,P)=\PP(\mathcal{
A}_1)/\PP(\mathcal{ B}).
\end{equation}

\begin{pf*}{Proof of Lemma~\ref{clustering-c}}
Here we consider the graph $G_s(n,m,\delta_x)$, where attribute sets of
vertices are of size $x$,
that is, $X_i=x$, $1\le i\le n$. We write ${\mathring{p}}=\PP(|D_{ij}|=s)$.

We note that (\ref{clustering-c-1}) follows from (\ref{ALPHA-H1}) and
the identities, which we prove below:
%
\begin{eqnarray}
\label{clusterlema-proof-1}  \PP(\mathcal{ B})&=&\PP(v_1\sim
v_2)\PP(v_1\sim v_3)={\mathring{p}}{}^2
\bigl(1+o(1) \bigr),
\\
\label{clusterlemma-proof-2} \PP(\mathcal{ A}_3)&\le&\PP(\mathcal{
A}_1)\le\PP(\mathcal{ A}_3)+o\bigl({
\mathring{p}}{}^2\bigr),
\\
\label{clusterlemma-proof-3}  \PP(\mathcal{ A}_3)&=&{\pmatrix{x
\cr
s}}^{-1}{\mathring{p}}{}^2 \bigl(1+o(1) \bigr).
\end{eqnarray}

We obtain the first identity of (\ref{clusterlema-proof-1}) from the
corresponding
identity of conditional probabilities
given $D_1$ that holds almost surely,
\begin{eqnarray*}
\PP(\mathcal{ B})=\PP(\mathcal{ B}|D_1)=\PP(v_1\sim
v_2|D_1)\PP(v_1\sim v_3|D_1)=
\PP(v_1\sim v_2)\PP(v_1\sim
v_3).
\end{eqnarray*}
The second identity of (\ref{clusterlema-proof-1}) follows from Lemma
\ref{sX1}. In order to show (\ref{clusterlemma-proof-3}), we write
$\PP(\mathcal{ A}_3)= {p_*}{\mathring{p}}$, where
${p_*}=\PP (D_3\cap(D_1\cup D_2)=D_{12}  |  |D_{12}|=s )$,
and evaluate
\begin{eqnarray*}
{p_*} = \frac{\bigl({m-(2x-s)\atop x-s}\bigr)}{\bigl({m\atop  x}\bigr)} = \frac{(m-2x+s)_{x-s}}{(m-x)_{x-s}}
\frac{{\mathring{p}}}{\bigl({x\atop s}\bigr)}=
\bigl(1+o(1) \bigr) \frac{{\mathring{p}}}{\bigl({x\atop s}\bigr)}.
\end{eqnarray*}
Let us prove (\ref{clusterlemma-proof-2}).
Since the event $\mathcal{ A}_3$ implies $\mathcal{ A}_2$, and $\mathcal{ A}_2$
implies $\mathcal{ A}_1$, we have
%
\begin{equation}
\label{22} \PP(\mathcal{ A}_1)= \PP(\mathcal{ A}_3)+\PP(
\mathcal{ A}_2\setminus\mathcal{ A}_3)+\PP (\mathcal{
A}_1\setminus\mathcal{ A}_2).
\end{equation}
It remains to show that
$\PP(\mathcal{ A}_2\setminus\mathcal{ A}_3), \PP(\mathcal{ A}_1\setminus\mathcal{
A}_2)=o({\mathring{p}}{}^2)$.
Noting that the event $\mathcal{ A}_1\setminus\mathcal{ A}_2$ implies at least
one of events
$\mathcal{ D}_{ij}$, $1\le i<j\le3$, we write, by the union bound and symmetry,
%
\begin{eqnarray}
\label{clusteringlemma-proof-4} \PP(\mathcal{ A}_1\setminus\mathcal{
A}_2)&\le&\sum_{1\le i<j\le3}\PP (\mathcal{
D}_{ij}\cap\mathcal{ A}_1 ) = 3\PP (\mathcal{
D}_{12}\cap\mathcal{ A}_1 )
\nonumber
\\[-8pt]
\\[-8pt]
\nonumber
&\le& 3\PP \bigl(\mathcal{
D}_{12}\cap\{v_1\sim v_3\} \bigr).
\end{eqnarray}
Now, conditioning on $D_1$, we obtain from inequalities of Lemma \ref
{sX1} that
\begin{eqnarray*}
\PP \bigl(\mathcal{ D}_{12}\cap\{v_1\sim
v_3\} | D_1 \bigr) & = & \PP (\mathcal{
D}_{12} | D_1 )\times\PP (v_1\sim
v_3 | D_1 )
\\
\nonumber
& \le & {\pmatrix{x
\cr
s+1}}^2{\pmatrix{m
\cr
s+1}}^{-1}\times {\pmatrix{x
\cr
s}}^2{\pmatrix{m
\cr
s}}^{-1}
\\
\nonumber
& = & \frac{(x-s)^2}{(m-s)(s+1)} {\pmatrix{x
\cr
s}}^4{\pmatrix{m
\cr
s}}^{-2}
\\
\nonumber
& = & o\bigl({\mathring{p}}{}^2\bigr). 
%
\end{eqnarray*}
We conclude that $\PP (\mathcal{ D}_{12}\cap\{v_1\sim v_3\}
)=o({\mathring{p}}{}^2)$.
Hence, $\PP(\mathcal{ A}_1\setminus\mathcal{ A}_2)=o({\mathring{p}}{}^2)$.
Next, observing that the event $\mathcal{ A}_2\setminus\mathcal{ A}_3$ implies
$\mathcal{ C}$, we write
%
\begin{equation}
\label{A1A2C} \PP(\mathcal{ A}_2\setminus\mathcal{ A}_3)
\le\PP(\mathcal{ C})
\end{equation}
and evaluate $\PP(\mathcal{ C})={\tilde p}_{1}{\tilde p}_{2}$, where
\begin{eqnarray*}
{\tilde p}_{1}=\PP\bigl(|D_{23}\setminus D_1|\ge1
| |D_{12}|=|D_{13}|=s\bigr),\qquad {\tilde p}_{2}=
\PP\bigl(|D_{12}|=|D_{13}|=s\bigr).
\end{eqnarray*}
Here ${\tilde p}_{1}$ is the probability that two independent subsets
$D_2\setminus D_1$ and $D_3\setminus D_1$ of $W\setminus D_1$
of size $|D_2\setminus D_1|=|D_3\setminus D_1|=x-s$ do intersect. By
Lemma~\ref{sX1},
${\tilde p}_{1}\le(x-s)^2(m-x)^{-1}$.
Hence, ${\tilde p}_{1}=o(1)$. Finally, the simple inequality
${\tilde p}_{2}\le\PP(\mathcal{ B})=(1+o(1)){\mathring{p}}{}^2$
implies $\PP(\mathcal{ C})=o({\mathring{p}}{}^2)$.
\end{pf*}

\begin{pf*}{Proof of Lemma~\ref{Lemma-clustering-1}} The proof goes
along the lines of the proof of Lemma~\ref{clustering-c}.

In view of (\ref{ALPHA-H1}), relation (\ref{clustering-c-2}) would
follow if we show that for any $\varepsilon\in(0, 1)$, we have, as
$m\to\infty$,
%
\begin{eqnarray}
\label{clustering-theorem-proof-1} \qquad\bigl(1-3\varepsilon^2+o(1)
\bigr)a_2a_1^2M^{-2}&\le&\PP(
\mathcal{ B})\le a_2a_1^2M^{-2},
\\
\label{clustering-theorem-proof-2} \PP(\mathcal{ A}_3)&\le&\PP(
\mathcal{ A}_1)\le\PP(\mathcal{ A}_3)+\bigl(
\varepsilon^2+o(1)\bigr) a_2a_1^2M^{-2},
\\
\label{clustering-theorem-proof-3} \bigl(1-2\varepsilon^2+o(1)
\bigr)a_1^3M^{-2}&\le&\PP(\mathcal{
A}_3) \le a_1^3M^{-2}.
\end{eqnarray}

We fix $\varepsilon$ and
introduce the indicator functions
${\mathbb I}_i=
{\mathbb I}_{\{X_i< \varepsilon\sqrt m\}}$ and write
${\overline{\mathbb I}}_i=
1-{\mathbb I}_i$
and
${\mathbb I}_*=
{\mathbb I}_1{\mathbb I}_2{\mathbb I}_3$.

Let us show (\ref{clustering-theorem-proof-1}).
We write ${\tilde\PP}(\mathcal{ B})={\tilde\PP}(v_1\sim v_2){\tilde\PP
}(v_1\sim v_3)$
and apply (\ref{sp}) to each probability ${\tilde\PP}(v_1\sim v_i)$,
$i=1,2$. We have
%
\begin{equation}
\label{PauliusDrun} Y_1^2Y_2Y_3M^{-2}(1-r){
\mathbb I}_* \le {\tilde\PP}(\mathcal{ B}){\mathbb I}_* \le {\tilde\PP}(\mathcal{
B}) \le Y_1^2Y_2Y_3M^{-2}.
\end{equation}
Here
\begin{eqnarray*}
r=1- \biggl(1-\frac{(X_1-s)_+(X_2-s)_+}{m-X_1+1} \biggr) \biggl(1-\frac
{(X_1-s)_+(X_3-s)_+}{m-X_1+1} \biggr).
\end{eqnarray*}
Observing that ${\mathbb I}_*=1$ implies $(X_i-s)_+\le\varepsilon
\sqrt{m}$ for each $i$, we write
${\mathbb I}_*r\le2\varepsilon^2(1-m^{-1/2})^{-1}$.
Now
(\ref{PauliusDrun}) implies the inequalities
%
\begin{eqnarray}
\label{B47}  Y_1^2Y_2Y_3M^{-2}
\bigl(1-2\varepsilon^2\bigl(1-m^{-1/2}\bigr)^{-1}-(1-{
\mathbb I}_*) \bigr) &\le& {\tilde\PP}(\mathcal{ B})
\nonumber
\\[-8pt]
\\[-8pt]
\nonumber
 &\le& Y_1^2Y_2Y_3M^{-2}.
\end{eqnarray}
Finally, taking the expected values and using the simple bounds
%
\begin{eqnarray}
\label{Rugpj23++} \E Y_1^2Y_2Y_3M^{-2}(1-{
\mathbb I}_*) &\le& \E Y_1^2Y_2Y_3M^{-2}({
\overline{\mathbb I}}_1+{\overline{\mathbb I}}_2+{
\overline{\mathbb I}}_3)
\nonumber
\\[-8pt]
\\[-8pt]
\nonumber
&=&o(1)a_2a_1^2M^{-2},
\end{eqnarray}
we obtain (\ref{clustering-theorem-proof-1}).
In the very last last step of (\ref{Rugpj23++})
we used (\ref{xyze}) in the form $\E Y_i^k{\overline{\mathbb I}}_i=o(a_k)$.

In order to show (\ref{clustering-theorem-proof-2}), we combine
(\ref{22}), (\ref{clusteringlemma-proof-4}) and (\ref{A1A2C})
with the
following inequalities:
%
\begin{eqnarray}
\label{Rugpj23+} \PP \bigl(\mathcal{ D}_{12}\cap\{v_1\sim
v_3\} \bigr)&\le& \bigl(\varepsilon^2+o(1) \bigr)
a_2a_1^2M^{-2},
\nonumber
\\[-8pt]
\\[-8pt]
\nonumber
 \PP(\mathcal{ C})
&\le& \bigl(\varepsilon^2+o(1) \bigr) a_2a_1^2M^{-2}.
\end{eqnarray}
Let us consider the first probability. Since the event $\mathcal{ D}_{12}$
implies $\{v_1\sim v_2\}$, we can write
%
\begin{equation}
\label{rugpj25-1} {\tilde\PP} \bigl(\mathcal{ D}_{12}\cap
\{v_1\sim v_3\} \bigr) \le {\tilde\PP} \bigl(\mathcal{
D}_{12}\cap\{v_1\sim v_3\} \bigr){\mathbb
I}_*+ {\tilde\PP} (\mathcal{ B} ) (1-{\mathbb I}_*).
\end{equation}
Next, we apply (\ref{sp}) to each probability of the right-hand side,
%
\begin{equation}
\label{F-L-namai} {\tilde\PP} \bigl(\mathcal{ D}_{12}\cap
\{v_1\sim v_3\} \bigr) \le \frac{\bigl({X_1\atop  s+1}\bigr)\bigl({X_2\atop  s+1}\bigr)}{\bigl({m\atop  s+1}\bigr)}
\frac{Y_1Y_3}{M},\qquad {\tilde\PP} (\mathcal{ B} )\le\frac{Y_1^2Y_2Y_3}{M^2},
\end{equation}
and write $\frac{\bigl({X_1\atop  s+1}\bigr)\bigl({X_2\atop  s+1}\bigr)}{\bigl({m\atop  s+1}\bigr)}=\frac{Y_1Y_2}{M}r'$, where
$r'=\frac{(X_1-s)_+(X_2-s)_+}{(m-s)(s+1)}$. Finally, collecting
inequalities~(\ref{F-L-namai}) and
$r'{\mathbb I}_*\le\varepsilon^2$
in
(\ref{rugpj25-1}) and then taking the expected value in~(\ref
{rugpj25-1}), we obtain the first bound of~(\ref{Rugpj23+}).
In order to prove the second bound of (\ref{Rugpj23+}), we write
%
\begin{equation}
\label{Rugpj23} {\tilde\PP}(\mathcal{ C})={\tilde p}_1{\tilde
p}_2\le{\tilde p}_1{\tilde p}_2{\mathbb
I}_*+{\tilde p}_2(1-{\mathbb I}_*),
\end{equation}
where ${\tilde p}_1={\tilde\PP} (|D_{23}\setminus D_1|\ge1
|  |D_{12}|=|D_{13}|=s )$ and
${\tilde p}_2={\tilde\PP}( |D_{12}|=|D_{13}|=s)$,
and apply (\ref{sp}) to ${\tilde p}_1$ and ${\tilde p}_2$,
\begin{eqnarray*}
{\tilde p}_1\le\frac{(X_2-s)_+(X_3-s)_+}{m-X_1},\qquad {\tilde p}_2\le
\frac{Y_1^2Y_2Y_3}{M^2}.
\end{eqnarray*}
We collect these inequalities in (\ref{Rugpj23}) and observe that
${\tilde p}_1{\mathbb I}_*\le\varepsilon^2/(1-m^{-1/2})$.
Now taking the expected value in (\ref{Rugpj23}), we obtain the second
bound of~(\ref{Rugpj23+}).


Let us show (\ref{clustering-theorem-proof-3}). We write $\mathcal{ A}_3$
in the form $\mathcal{ L}_1\cap\mathcal{ L}_2$, where
$\mathcal{ L}_1=\{D_3\cap(D_1\cup D_2)=D_{12}\}$, $\mathcal{ L}_2=\{|D_{12}|=s\}$,
and decompose ${\tilde\PP}(\mathcal{ A}_3)={\tilde p}_3{\tilde p}_4$. Here
${\tilde p}_3={\tilde\PP}(\mathcal{ L}_1|\mathcal{ L}_2)$
and
${\tilde p}_4={\tilde\PP}(\mathcal{ L}_2)$ satisfy
%
\begin{eqnarray}
\label{up2}  (1-{\tilde r}_3)Y_3M^{-1}&\le&{
\tilde p}_3\le Y_3M^{-1}, \qquad {\tilde
r}_3:=\frac{(X_1+X_2-2s)(X_3-s)}{m-X_3},
\\
\label{up3}\qquad\quad  (1-{\tilde r}_4)Y_1Y_2M^{-1}
&\le&{\tilde p}_4\le Y_1Y_2M^{-1},\qquad {
\tilde r}_4:=\frac{(X_1-s)_+(X_2-s)_+}{m-X_1+1}.
\end{eqnarray}
Note that (\ref{up3}) is an immediate consequence of (\ref{sp}). In
order to get (\ref{up2}), we write
%
\begin{equation}
\label{up1} {\tilde p}_3= \frac{\bigl({m-X_1-X_2+s\atop  X_3-s}\bigr)}{\bigl({m\atop  X_3}\bigr)} = \theta
\frac{Y_3}{M}, \qquad\theta:=\frac{(m-X_1-X_2+s)_{X_3-s}}{(m-s)_{X_3-s}},
\end{equation}
and apply to $\theta$ the chain of inequalities
\begin{eqnarray*}
\frac{(a)_r}{(a+b)_r}>\frac{(a-r)^r}{(a+b-r)^r}= \biggl(1-\frac
{b}{a+b-r}
\biggr)^r\ge1-\frac{br}{a+b-r}.
\end{eqnarray*}

Combining (\ref{up1}), (\ref{up2}), (\ref{up3}) and
${\tilde\PP}(\mathcal{ A}_3)={\tilde p}_3{\tilde p}_4$, we write
%
\begin{equation}
\label{rugpj25+2}\qquad (1-{\tilde r}_3) (1-{\tilde r}_4)Y_1Y_2Y_3M^{-2}{
\mathbb I}_* \le {\tilde\PP}(\mathcal{ A}_3){\mathbb I}_* \le {\tilde
\PP}(\mathcal{ A}_3)\le Y_1Y_2Y_3M^{-2}.
\end{equation}
Finally, invoking in (\ref{rugpj25+2}) the inequalities
${\tilde r}_i{\mathbb I}_*\le\varepsilon^2(1-m^{-1/2})^{-1}$ and then
taking the expected values,
we obtain (\ref{clustering-theorem-proof-3}).
\end{pf*}

\begin{pf*}{Proof of Theorem~\ref{Theorem-clustering-1+}}
In the case where (ii$'$) holds, we apply Lemma~\ref{Lemma-clustering-1}.
Note that conditions (i) and (ii$'$) imply the uniform integrability of
the sequences of random variables $\{Z^2_{m1}\}_m$ and
$\{Z_{m1}\}_m$; that is, for $k=1,2$, we have
%
\begin{equation}
\label{uniform-integrability} \forall\varepsilon>0, \exists \Delta>0 {\mbox{ such
that }} \forall m\qquad \E Z^k_{m1}{\mathbb I}_{\{Z_{m1}>\Delta\}}<
\varepsilon;
\end{equation}
see Lemma~\ref{auxiliarylema}.
Hence, (\ref{xyze}) holds, and Theorem~\ref{Theorem-clustering-1+}
follows from
Lemma~\ref{Lemma-clustering-1}.

Now assume that $\E Z^2=+\infty$. In order to prove that $\alpha=o(1)$, we
show [see~(\ref{ALPHA-H1})] that
%
\begin{equation}
\label{2011-12AB} \PP(\mathcal{ A}_1)=O\bigl(n^{-2}\bigr)
\quad{\mbox{and}}\quad \liminf_mn^{2}\PP(\mathcal{ B})=+\infty.
\end{equation}
Before the proof of (\ref{2011-12AB}), we introduce some notation. Denote
\begin{eqnarray*}
h\to\varphi(h)=\sup_m \E X_{m1}^s{\mathbb
I}_{\{X_{m1}>h\}},\qquad \varepsilon_h=\varphi\bigl(h^{1/2}
\bigr),\qquad \delta^{2s}_h=\max\bigl\{h^{-1},
\varepsilon_h\bigr\},
\end{eqnarray*}
and observe, that (i), (ii) and (\ref{beta+++}) imply $\varphi(h)=o(1)$
as $h\to\infty$.
Hence, $\delta_h=o(1)$ as $h\to+\infty$. Furthermore, by Chebyshev's inequality,
%
\begin{eqnarray}
\label{H+H} \PP(X_{m1}>m\delta_m) &\le& (m
\delta_m)^{-s}\varphi(m\delta_m) \le (m
\delta_m)^{-s}\varepsilon_m \le m^{-s}
\varepsilon_m^{1/2}
\nonumber
\\[-8pt]
\\[-8pt]
\nonumber
& =& o\bigl(m^{-s}\bigr).
\end{eqnarray}
For $k=2,3$ and $t=1,3$, denote $D_k^*:=D_k\setminus D_1$ and
introduce the events $\mathcal{ H}_t=\{|D_t|\le m\delta_m\}$,
\begin{eqnarray*}
 \mathcal{ Q}_i&=&\bigl\{|D_1\cap D_2
\cap D_3|=i\bigr\},\qquad \mathcal{ Q}_i^*=\bigl\{\bigl|D_2^*
\cap D_3^*\bigr|\ge s-i\bigr\},\qquad 0\le i\le s-1,
\\
 Q_s&=&\bigl\{|D_1\cap D_2\cap
D_3|\ge s\bigr\},\qquad \mathcal{ Q}_s^*=\bigl\{\bigl|D_2^*
\cap D_3^*\bigr|\ge0\bigr\}, \qquad\mathcal{ H}=\mathcal{ H}_1\cap
\mathcal{ H}_3.
\end{eqnarray*}
Let us prove the first bound of (\ref{2011-12AB}). We denote $R=\PP
(\mathcal{ A}_1\cap{\overline\mathcal{ H}})$ and write
%
\begin{equation}
\label{A-+-H} \PP(\mathcal{ A}_1) = \PP(\mathcal{ A}_1
\cap\mathcal{ H})+R={\tilde p}_0+\cdots+{\tilde p}_s+R,
\end{equation}
where ${\tilde p}_i=\PP(\mathcal{ B}\cap\mathcal{ Q}_i\cap\mathcal{ Q}_i^*\cap
\mathcal{ H})$. Note that the remainder $R$ is negligibly small,
\begin{eqnarray*}
R 
\le\PP(v_2\sim v_3)\PP({\overline\mathcal{
H}_1}) + \PP(v_1\sim v_2)\PP({\overline
\mathcal{ H}_3})=o\bigl(m^{-2s}\bigr).
\end{eqnarray*}
In the last step, we used (\ref{H+H}) and the inequality $\PP(v_2\sim
v_3)\le a_1^2M^{-1}$, which follows from~(\ref{sp}).
Let us estimate
${\tilde p}_{i}
=
\PP(\mathcal{ Q}_i^*|\mathcal{ B}\cap\mathcal{ Q}_i\cap\mathcal{ H})\PP(\mathcal{ B}\cap
\mathcal{ Q}_i\cap\mathcal{ H})
=:
{\tilde p}_{i1}{\tilde p}_{i2}$.
Note that, on the event $\mathcal{ B}\cap\mathcal{ H}_1$, the sets
$D_2^*$ and $D_3^*$ are random subsets of $W\setminus D_1$ of sizes
$|D_k^*|\le X_k-s$, where
$|W\setminus D_1|\ge m(1-\delta_m)=:m'$ is of order $(1-o(1))m$. By
(\ref{sp}),
the probability that their intersection has at least $s-i$ elements
is bounded from above by $\E\bigl({X_2\atop  s-i}\bigr)\bigl({X_3\atop  s-i}\bigr){\bigl({m'\atop  s-i}\bigr)}^{-1}$.
Hence, ${\tilde p}_{i1}=O(m^{i-s})$.
Next, we estimate
\begin{eqnarray*}
{\tilde p}_{i2} &\le& \PP\bigl(\{v_1\sim v_2\}
\cap\mathcal{ Q}_i\cap\mathcal{ H}_3\bigr) = \PP
\bigl(\bigl|D_2^{**}\cap D_1^{**}\bigr|\ge s-i
|\mathcal{ Q}_i\cap\mathcal{ H}_3 \bigr)\PP(\mathcal{
Q}_i\cap\mathcal{ H}_3)\\
& =:&{\tilde p}_{i2}'
{\tilde p}_{i2}''.
\end{eqnarray*}
Here the random subsets
$D_k^{**}:=D_k\setminus D_3$ of $W\setminus D_3$ are of sizes at most
$X_k-i\le X_k$, $k=1,2$.
Since on the event $\mathcal{ H}_3$ we have $|W\setminus D_3|\ge m'$,
inequality (\ref{sp}) implies
${\tilde p}_{i2}'\le\E\bigl({X_1\atop  s-i}\bigr)\bigl({X_2\atop  s-i}\bigr){\bigl({m'\atop s-i}\bigr)}^{-1}=O(m^{i-s})$.
Finally, we estimate ${\tilde p}_{i2}''$. Let $S$ count the number of
subsets $A_i$ of $W$ of size $i$
covered
by the intersection $D_1\cap D_2\cap D_3$, that is,
$S=\sum_{A_i\subset W}{\mathbb I}_{\{A_i\subset D_1\cap D_2\cap D_3\}}$. We have
\begin{eqnarray*}
{\tilde p}_{i2}'' \le \PP(S\ge1)\le\E S =
\pmatrix{m
\cr
i} \bigl(\PP(A_i\in D_1)
\bigr)^3={\pmatrix{m
\cr
i}}^{-2} \biggl(\E{\pmatrix{X_1
\cr
i}} \biggr)^3.
\end{eqnarray*}
Hence, ${\tilde p}_{i2}''=O(m^{-2i})$. Collecting these bounds, we
conclude that
\begin{eqnarray*}
{\tilde p}_{i} = {\tilde p}_{i1}{\tilde
p}_{i2} \le {\tilde p}_{i1}{\tilde p}_{i2}'{
\tilde p}_{i2}'' = O\bigl(m^{-2s}
\bigr).
\end{eqnarray*}
Now the first bound of (\ref{2011-12AB}) follows from (\ref{A-+-H}).

Let us prove the second identity of (\ref{2011-12AB}).
Given $t<m$, denote ${\mathbb I}_{it}={\mathbb I}_{\{X_i\le t\}}$
and write [see (\ref{sp})]
%
\begin{eqnarray}
\label{BTTB} \PP(\mathcal{ B})& =& \E{\tilde\PP}(\mathcal{ B}) \ge \E{\mathbb
I}_{1t}{\mathbb I}_{2t}{\mathbb I}_{3t}{\tilde
\PP}(\mathcal{ B})
\nonumber
\\[-8pt]
\\[-8pt]
\nonumber
&\ge& M^{-2}\E{\mathbb I}_{1t}{\mathbb
I}_{2t}{\mathbb I}_{3t}Y_1^2Y_2Y_3
\bigl(1-t^2/(m-t)\bigr)^2.
\end{eqnarray}
We have, by (i) and (ii), that
$\E{\mathbb I}_{it}Y_i\to a_1$ as $m,t\to+\infty$. Furthermore,
(i) and $\E Z^2=\infty$ implies $\E Y_1^2{\mathbb I}_{1t}\to+\infty$ as
$m,t\to\infty$.
Therefore, letting $t=t_m\to+\infty$ and $t_m=o(m^{1/2})$, we obtain
(\ref{2011-12AB}).
\end{pf*}

\begin{pf*}{Proof of Remark~\ref{rem2}}
Here we assume that $\beta_m:=\bigl({m\atop s}\bigr)n^{-1}$ tend to $+\infty$ as
$m,n\to\infty$. In the case where
(i), (ii) and (ii$'$) hold, the relation $\alpha_s(m,P)=o(1)$ is an
immediate consequence of Lemma
\ref{Lemma-clustering-1}.
In the case where (i) and (ii) hold and $\E Z^2=\infty$, we show
(\ref{2011-12AB})
by the same argument as that used in the proof of Theorem \ref
{Theorem-clustering-1+} above. In particular,
in the proof of the first bound of (\ref{2011-12AB}), we choose
$\delta_m=\beta_m^{-1/(4s)}$, and, for $\tau_m=\bigl({m\delta_m\atop  s}\bigr)\beta_m^{-1/2}$, we
write [cf. (\ref{H+H})]
\begin{eqnarray*}
\PP(X_{m1}>m\delta_m)=\PP(Z_{m1}>
\tau_m)\le\tau_m^{-1}\E Z_{m1}{
\mathbb I}_{\{Z_{m1}>\tau_m\}}=o\bigl(\tau_m^{-1}\bigr).
\end{eqnarray*}
Here $\tau_m^{-1}\le c(s) m^{-s/2}n^{-1/2}\delta_m^{-s}=o(n^{-1})$.
In the proof of the second identity of~(\ref{2011-12AB}), we choose $t=t_m=m^{1/2}n^{-1/(4s)}$ in (\ref
{BTTB}). Finally, from (\ref{2011-12AB}) and~(\ref{ALPHA-H1}) we deduce the relation
$\alpha_s(m,P)=o(1)$.
\end{pf*}

\begin{pf*}{Proof of Theorem~\ref{Theorem-clustering-degree}} We
start with some notation.
By $f_i(\lambda)=e^{-\lambda}\lambda^i/i!$ we denote Poisson's
probability. Observing that the absolute value
of the derivative
of the function $\lambda\to f_i(\lambda)$ is bounded by $1$,
we write
$|f_i(\lambda_1)-f_i(\lambda_2)|\le|\lambda_1-\lambda_2|$, by
the mean value theorem. We denote
\begin{eqnarray*}
H_{1}=M^{-2}Y_1Y_2Y_3
=n^{-2}\beta_m^{-1/2}Z_1Z_2Z_3,\qquad
H_2=n^{-2}Z_1^2Z_2Z_3,
\end{eqnarray*}
$\mu_1=\E Z_1$, and $\beta_m=M/n$, and introduce the event
$\mathcal{ K}= \{\sum_{j=4}^{n}{\mathbb I}_{\{v_1\sim v_j\}}=k-2 \}
$. We also note that, for $i=1,2$,
the moments $\E Z_1^i$ and $a_i=\E Y_1^i=\beta_m^{-i/2}\E Z_1^i$ are
bounded from above and that they are bounded
away from zero as $m,n\to\infty$.

Let us prove (\ref{clustering-degree}). Given integer $k\ge2$, we write
%
\begin{equation}
\label{alpha-cluster-degree} \alpha^{[k]} = \frac{\PP(\mathcal{ A}_1\cap\{ d(v_1)=k \} )} {\PP( \mathcal{ B}\cap\{d(v_1)=k\})} =
\frac{p^{*}}{p^{**}}, 
\end{equation}
where $p^{*}=\PP(\mathcal{ A}_1\cap\mathcal{ K})$ and $p^{**}=\PP(\mathcal{ B}\cap
\mathcal{ K})$.
We shall show that, for any $\varepsilon\in(0,1)$,
%
\begin{eqnarray}
\label{r-1}  p^*&=&n^{-2}\beta_m^{-1/2}
\mu_1(k-1)\E f_{k-1}(\mu_1
Z_1)+R_1,
\\
\label{r-2}  p^{**} &=& n^{-2}(k-1)k\E f_k(
\mu_1Z_1) +R_2,
\end{eqnarray}
where the remainder terms are negligibly small,
%
\begin{eqnarray*}
|R_1|&\le& n^{-2}c(\beta) \bigl(a_2a_1^2+a_1^2+a_1^3
\bigr) \bigl(\varepsilon +o(1) \bigr),\\
 |R_2|&\le& n^{-2}c(
\beta) \bigl(a_2a_1^2+a_1^2
\bigr) \bigl(\varepsilon+o(1) \bigr).
\end{eqnarray*}
Here the number $c(\beta)$ depends only on $\beta$. Note that (i) and
(ii) imply
$\mu_1\to\mu$ and the convergence in distribution $\mu_1Z_1\to\mu Z$. Hence,
we have the convergence of the expectations
$\E f_i(\mu_1Z_1)\to\E f_i(\mu Z)=p_i$. Now (\ref{r-1}), (\ref{r-2})
and (\ref{beta+++}) imply
$p^*=n^{-2}\beta^{-1/2}\mu(k-1)p_{k-1}+o(n^{-2})$ and $p^{**} =
n^{-2}(k-1)kp_k+o(n^{-2})$.
Substitution of these identities into
(\ref{alpha-cluster-degree}) shows (\ref{clustering-degree}).
%

In the remaining part of the proof, we show (\ref{r-1}) and~(\ref{r-2}) for a given $\varepsilon\in(0,1)$.
For this purpose, we write
%
\begin{equation}
\label{p*p**} p^*=\E {\tilde\PP}(\mathcal{ A}_1){\tilde\PP}(
\mathcal{ K}), \qquad p^{**}=\E {\tilde\PP}(\mathcal{ B}) {\tilde\PP}(\mathcal{
K})
\end{equation}
and show that the conditional probabilities
${\tilde\PP}(\mathcal{ A}_1)$, ${\tilde\PP}(\mathcal{ B})$ and ${\tilde\PP
}(\mathcal{ K})$
can be approximated by
%
\begin{equation}
\label{approximate-p} {\tilde\PP}(\mathcal{ K})\approx f_{k-2}(
\mu_1 Z_1),\qquad {\tilde\PP}(\mathcal{ A}_1)
\approx H_1,\qquad {\tilde\PP}(\mathcal{ B})\approx H_2.
\end{equation}
Note that substitution
of (\ref{approximate-p}) into (\ref{p*p**}) gives the leading terms of
(\ref{r-1}) and~(\ref{r-2}).
In the remaining part of the proof, we show the validity of such an
approximation.

Note that since conditions of Theorem~\ref{Theorem-clustering-degree}
are more restrictive than
those of Lemma~\ref{Lemma-clustering-1} (see the proof of Theorem \ref
{Theorem-clustering-1+}), we can rightfully use the inequalities
established in the proof of Lemma~\ref{Lemma-clustering-1}.

\textit{Approximation of} ${\tilde\PP}(\mathcal{ K})$. We first establish an
upper bound for
$\Delta={\tilde\PP}(\mathcal{ K})-f_{k-2}(\mu_1 Z_1)$,
%
\begin{equation}
\label{Delta-1} |\Delta|\le2{\tilde\kappa}_1+{\tilde
\kappa}_2+|{\tilde\mu}-\mu_1|Z_1.
\end{equation}
Here
\begin{eqnarray*}
{\tilde{\mu}}&=&n^{-1}\sum_{i=4}^nZ_i,\qquad
{\tilde{\kappa}}_1=n^{-2}\sum
_{i=4}^n Z_1^2
Z_i^2,\\
 {\tilde{\kappa}}_2&=&n^{-1}
\sum_{i=4}^n Z_1
Z_i \frac{(X_1-s)(X_i-s)} {m-X_1+1}.
\end{eqnarray*}
We show (\ref{Delta-1}) in two steps.
In the first step, we apply Le Cam's inequality (\ref{LeCam}) to the sum
$\sum_{i=4}^{n}{\mathbb I}_{\{v_1\sim v_i\}}$ of conditionally
independent\vadjust{\goodbreak} random indicators,
%
\begin{equation}
\label{r13-7} \Biggl|{\tilde\PP}(\mathcal{ K})-f_{k-2} \Biggl(\sum
_{i=4}^n{\tilde\PP}(\mathcal{ E}_{1i})
\Biggr) \Biggr| \le 2\sum_{i=4}^n \bigl({\tilde
\PP}(\mathcal{ E}_{1i}) \bigr)^2.
\end{equation}
%
In the second step, using the mean value theorem, we replace
the argument of $f_{k-2}$ by $\mu_1 Z_1$ in~(\ref{r13-7}). The
approximation error estimate (\ref{Delta-1}) now
follows from the inequalities
%
\begin{equation}
\label{r13-8} \sum_{i=4}^n \bigl({\tilde
\PP}(\mathcal{ E}_{1i}) \bigr)^2\le{\tilde
\kappa}_1,\qquad {\tilde\mu}Z_1-{\tilde\kappa}_2
\le \sum_{i=4}^n{\tilde\PP}(\mathcal{
E}_{1i})\le{\tilde\mu}Z_1,
\end{equation}
which are simple consequences of the inequalities of Lemma~\ref{sX1}.

Next, using (\ref{Delta-1}), we replace ${\tilde P}(\mathcal{ K})$ by
$f_{k-2}(\mu_1 Z_1)$ in (\ref{p*p**})
and show that the error of the
replacement is negligibly small. We denote
${\mathbb I}_{A}={\mathbb I}_{\{Z_{1}\le A\}}$ and ${\overline{\mathbb
I}}_{A}=1-{\mathbb I}_{A}$ and observe, that
in view of~(\ref{uniform-integrability}), we can find
$A>1$ such that, uniformly in $m$,
%
\begin{equation}
\label{choice-of-A} \E Z_{1}^k{\overline{\mathbb
I}}_{A}<\varepsilon\E Z_{1}^k,\qquad k=1,2.
\end{equation}
We write
%
\begin{eqnarray}
\label{K-A}  {\tilde\PP}(\mathcal{ A}_1){\tilde\PP}(\mathcal{ K})
&=& {\tilde\PP}(\mathcal{ A}_1)f_{k-2}(\mu_1
Z_1){\mathbb I}_A+r_1+r_2,
\\
\label{K-B}  {\tilde\PP}(\mathcal{ B}){\tilde\PP}(\mathcal{ K}) &= &{\tilde\PP}(
\mathcal{ B})f_{k-2}(\mu_1 Z_1){\mathbb
I}_A+r_3+r_4,
\end{eqnarray}
where
\begin{eqnarray*}
r_1&=& {\tilde\PP}(\mathcal{ A}_1){\tilde\PP}(\mathcal{
K}){\overline{\mathbb I}}_A,\qquad r_2={\tilde\PP}(\mathcal{
A}_1)\Delta{\mathbb I}_A,\\
r_3&=&{\tilde\PP}(
\mathcal{ B}){\tilde\PP}(\mathcal{ K}){\overline{\mathbb I}}_A,\qquad
r_4={\tilde\PP}(\mathcal{ B})\Delta{\mathbb I}_A,
\end{eqnarray*}
and show that
%
\begin{eqnarray}
\label{rrrr+} \E r_i&\le&\varepsilon M^{-2}
a_2a_1^2 \qquad{\mbox{for }} i=1,3,
\nonumber
\\[-8pt]
\\[-8pt]
\nonumber
\E|r_j|&=&o\bigl(M^{-1}n^{-1}\bigr)a_1^2
\qquad{\mbox{for }} j=2,4.
\end{eqnarray}
To prove (\ref{rrrr+}), we write, using the inequalities
${\tilde\PP}(\mathcal{ A}_1)\le{\tilde\PP}(\mathcal{ B})$ and
${\tilde\PP}(\mathcal{ B})\le H_2$ [see (\ref{sp})],
%
\begin{equation}
\label{rrrr} r_i\le H_2{\overline{\mathbb
I}}_A \qquad{\mbox{for }} i=1,3 \quad{\mbox{and}} \quad |r_j|\le
H_2|\Delta|{\mathbb I}_A \qquad{\mbox{for }} j=2,4
\end{equation}
and notice that, for $i=1,3$, inequalities (\ref{rrrr+}) follow from
(\ref{choice-of-A}) and,
for $j=2,4$,
inequalities (\ref{rrrr+}) follow from the bound $\E Z_1^2|\Delta|
{\mathbb I}_A=o(1)$.
Let us prove this bound.
We note that (\ref{Delta-1}), combined with the simple inequality
$|\Delta|\le1$, implies
\[
|\Delta|\le1\wedge(2{\tilde\kappa}_1)+1\wedge\bigl(Z_1|{
\tilde\mu}-\mu_1|\bigr)+1\wedge{\tilde\kappa}_2.
\]
Here we use the notation $a\wedge b=\min\{a,b\}$. Hence, for $A>1$, we
can write
%
\begin{equation}
\label{Delta-3} Z_1^2|\Delta|{\mathbb I}_A
\le c_A \bigl( 1\wedge\kappa^*_1+1\wedge|{\tilde\mu}-
\mu_1|+1\wedge\kappa^*_2 \bigr),\vadjust{\goodbreak}
\end{equation}
where the number $c_A>0$ does not depend on $m$, and where
\begin{eqnarray*}
\kappa^*_1:= n^{-2}\sum_{i=4}^nZ_i^2,\qquad
\kappa_2^*:=n^{-1}m^{-1}\sum
_{i=4}^n Z_i (X_i-s).
\end{eqnarray*}
The bound $\E Z_1^2|\Delta| {\mathbb I}_A=o(1)$ follows from (\ref
{Delta-3}) and from the bounds
%
\begin{equation}
\label{Delta-41} \kappa^*_1=o_P(1),\qquad {\tilde\mu}-
\mu_1=o_P(1), \qquad\kappa_2^*=o_P(1).
\end{equation}
The first and second bounds of (\ref{Delta-41}) are shown in Lemma \ref
{auxiliarylema}.
To obtain the third bound, we write $\E\kappa_2^*\le m^{-1}\E
Z_4X_4=o(1)$; see (i), (ii$'$).

\textit{Approximation of} ${\tilde\PP}(\mathcal{A}_1)$.
Now we replace ${\tilde\PP}(\mathcal{ A}_1)$ by $H_{1}$ in the right-hand
side of (\ref{K-A})
and show that the error of the replacement is negligibly small.
We proceed in two steps. First, we approximate
${\tilde\PP}(\mathcal{ A}_1)\approx{\tilde\PP}(\mathcal{ A}_3)$,
and, second, we approximate ${\tilde\PP}(\mathcal{ A}_3)\approx H_{1}$.
In the first step, we combine the inequalities
[see (\ref{clusteringlemma-proof-4}), (\ref{A1A2C}) and~(\ref{Rugpj23+})]
\begin{eqnarray*}
{\tilde\PP}(\mathcal{ A}_3) &\le& {\tilde\PP}(\mathcal{
A}_1) = {\tilde\PP}(\mathcal{ A}_3)+{\tilde\PP}(
\mathcal{ A}_1\setminus\mathcal{ A}_2)+{\tilde\PP}(
\mathcal{ A}_2\setminus\mathcal{ A}_3),
\\
 \E{\tilde\PP}(\mathcal{ A}_1\setminus\mathcal{
A}_2) &\le &3\PP \bigl(\mathcal{ D}_{12}\cap
\{v_1\sim v_3\} \bigr) \le 3 \bigl(\varepsilon^2+o(1)
\bigr)M^{-2}a_2a_1^2,
\\
 \E{\tilde\PP}(\mathcal{ A}_2\setminus\mathcal{
A}_3) &\le &\PP(\mathcal{ C}) \le \bigl(\varepsilon^2+o(1)
\bigr)M^{-2}a_2a_1^2
\end{eqnarray*}
and show that
%
\begin{eqnarray}
\label{r-6} {\tilde\PP}(\mathcal{ A}_3)f_{k-2}(
\mu_1 Z_1){\mathbb I}_A &\le& {\tilde\PP}(
\mathcal{ A}_1)f_{k-2}(\mu_1
Z_1){\mathbb I}_A
\nonumber
\\[-8pt]
\\[-8pt]
\nonumber
&\le&{\tilde\PP}(\mathcal{
A}_3)f_{k-2}(\mu_1 Z_1){\mathbb
I}_A +r_5,
\end{eqnarray}
where $\E r_5\le4 (\varepsilon^2+o(1) )M^{-2}a_2a_1^2$.
In the second step, we apply (\ref{rugpj25+2}) and write
%
\begin{equation}
\label{r-7} \qquad H_{1}f_{k-2}(\mu_1
Z_1){\mathbb I}_A- r_6 \le {\tilde\PP}(
\mathcal{ A}_3)f_{k-2}(\mu_1 Z_1){
\mathbb I}_A
\le H_{1}f_{k-2}(
\mu_1 Z_1){\mathbb I}_A,
\end{equation}
where $\E r_6\le (2\varepsilon^2+o(1) )M^{-2}a_1^3$.
Finally, using the inequality $1-{\overline{\mathbb I}}_A ={\mathbb
I}_A\le1$, we write
%
\begin{equation}
\label{H-A} \E H_{1}f_{k-2}(\mu_1
Z_1)-r_7 = \E H_{1}f_{k-2}(
\mu_1 Z_1){\mathbb I}_A\le\E
H_{1}f_{k-2}(\mu_1 Z_1),
\end{equation}
where
$r_7=\E H_{1}f_{k-2}(\mu_1 Z_1){\overline{\mathbb I}}_A
\le\E H_{1}{\overline{\mathbb I}}_A\le \varepsilon M^{-2}a_1^3$.
In the last step, we\break used~(\ref{choice-of-A}).
Finally, from (\ref{p*p**}), (\ref{K-A}), (\ref{r-6}), (\ref{r-7}) and (\ref{H-A}) we
obtain~(\ref{r-1}).

\textit{Approximation of} ${\tilde\PP}(\mathcal{ B})$.
Now we replace
${\tilde\PP}(\mathcal{ B})$ by $H_{2}$ in the right-hand side of (\ref{K-B}).
For this purpose, we combine (\ref{PauliusDrun}) and (\ref{B47})
with the inequality $f_{k-2}(\mu Z_1)\le1$ and obtain
%
\begin{equation}
\label{r-5} \qquad H_{2}f_{k-2}(\mu_1
Z_1){\mathbb I}_A-r_8 \le {\tilde\PP}(
\mathcal{ B})f_{k-2}(\mu_1 Z_1){\mathbb
I}_A \le H_{2}f_{k-2}(\mu_1
Z_1){\mathbb I}_A,
\end{equation}
where $r_8\le (3\varepsilon+(1-{\mathbb I}_{*}) )H_2$ is
negligibly small, because
$\E r_8\le(3\varepsilon+o(1))\times  M^{-2}a_2a_1^2$. Next, we proceed as in
(\ref{H-A}) above and write
%
\begin{equation}\qquad
\label{H-B} \E H_{2}f_{k-2}(\mu_1
Z_1)-r_9= \E H_{2}f_{k-2}(
\mu_1 Z_1){\mathbb I}_A\le\E
H_{2}f_{k-2}(\mu_1 Z_1),
\end{equation}
where $r_9=\E H_{2}f_{k-2}(\mu_1 Z_1){\overline{\mathbb I}}_A
\le\E H_{2}{\overline{\mathbb I}}_A\le \varepsilon M^{-2}a_2a_1^2$.
Finally, from (\ref{p*p**}), (\ref{K-B}), (\ref{r-5}) and (\ref{H-B}) we obtain (\ref{r-2}).
\end{pf*}

\subsection{Passive graph}\label{sec4.3}
\mbox{}
\begin{pf*}{Proof of Theorem~\ref{PTh-0}}
%
%
%
The proof consists of two parts. In the first part, we
establish the convergence
of the Fourier transforms
%
\begin{equation}
\label{f-d*} \lim_{m\to\infty}\E e^{itL}=\E e^{itd_*}.
\end{equation}

In the second part, we
show that, as $m\to\infty$,
%
\begin{equation}
\label{d=L} \PP(d\not=L)=o(1).
\end{equation}

\textit{Part} 1. Here we prove (\ref{f-d*}).
%
Before the proof, we introduce some notation.
We denote ${\mathbb I}_j(w)={\mathbb I}_{\{w\in D_j\}}$ and write
%
\begin{equation}
\label{L=L} L=(X_1-1)_+{\mathbb I}_1(w_1)+
\cdots+(X_n-1)_+{\mathbb I}_n(w_1).
\end{equation}
Here we use the notation $a_+=\max\{0,a\}$. Given ${\mathbb
X}=(X_1,\ldots, X_n)$, we generate
independent Poisson random variables $\eta_1(X_1),\ldots, \eta_n(X_n)$
with mean values
$\E(\eta_j(X_j)|X_j)=m^{-1}X_j$. We also generate independent random variables
$Y_1({\mathbb X}), Y_2({\mathbb X}),\ldots$ with the common probability
distribution
defined as follows. In the case where $S(X)=X_1+\cdots+X_n$ is positive,
we put
\begin{eqnarray*}
\PP \bigl(Y_1({\mathbb X})=j |{\mathbb X} \bigr)=(j+1)
\frac{n_{j+1}}{S(X)},\qquad n_{j+1}=\sum_{k=1}^n{
\mathbb I}_{\{X_k=j+1\}},\qquad j=0,1,\ldots.
\end{eqnarray*}
For $S(X)=0$, we put $Y_1({\mathbb X})\equiv0$.
By $\eta_1,\ldots, \eta_n$ (resp., $Y_1,Y_2,\ldots$)
we denote the corresponding unconditional random
variables, that is, $\eta_j$ is the outcome of a two-step procedure: first,
we generate ${\mathbb X}$ and then
generate $\eta_j(X_j)$.
Let $\eta({\mathbb X})=\eta_1(X_1)+\cdots+\eta_n(X_n)$ and $\eta=\eta_1+\cdots+\eta_n$.
Introduce the random variables
%
\begin{equation}
\label{tau-xi} \tau=(X_1-1)_+\eta_1+
\cdots+(X_n-1)_+\eta_n,\qquad \xi=\sum
_{j=1}^\eta Y_j
\end{equation}
and observe that, given ${\mathbb X}$, their conditional distributions
$P_{\tau}({\mathbb X})$ and $P_{\xi}({\mathbb X})$ are compound Poisson
and that $P_{\tau}({\mathbb X})\equiv P_{\xi}({\mathbb X})$.
Here and below in this proof, by $P_\zeta({\mathbb X})$ (and $\E_{\mathbb X}\zeta$) we denote the
conditional distribution (and conditional expectation) of a random
variable $\zeta$ given
${\mathbb X}$ .

Let us prove (\ref{f-d*}). For this purpose, we write $\E e^{itL}-\E
e^{itd_*}= I_1+I_2$, where\vadjust{\goodbreak}
$I_1=\E e^{itL}-\E e^{it\xi}$ and $I_2=\E e^{it\xi}-\E e^{itd_*}$, and
show that
$I_1, I_2=o(1)$.
In order to show the first bound, we write
\begin{eqnarray*}
I_1= \E \bigl(\E_{\mathbb X}e^{itL}-
\E_{\mathbb X}e^{it\xi} \bigr) =\E{\tilde\Delta},\qquad {\tilde\Delta}=
\bigl(\E_{\mathbb X}e^{itL}-\E_{\mathbb X}e^{it\tau
} \bigr)
\end{eqnarray*}
and invoke the inequalities
%
\begin{equation}
\label{TV} 2^{-1}\bigl\llvert \E_{\mathbb X}e^{itL}-
\E_{\mathbb X}e^{it\tau}\bigr\rrvert \le d_{\mathrm{TV}}
\bigl(P_L({\mathbb X}),P_{\tau}({\mathbb X}) \bigr)
\le2m^{-2}\sum_{j=1}^nX_j^2.
\end{equation}
Here $d_{\mathrm{TV}}$ denotes the total variation distance.
Indeed, we obtain from (\ref{TV}) that, for any $\varepsilon>0$,
\begin{eqnarray*}
|I_1|\le\varepsilon+2\PP\bigl(|{
\tilde\Delta}|>\varepsilon\bigr)\le\varepsilon +2\PP \Biggl(4m^{-2}\sum
_{j=1}^nX_j^2>
\varepsilon \Biggr)=\varepsilon+o(1).
\end{eqnarray*}
In the last step, we applied (\ref{ThirdFact}). Hence, $I_1=o(1)$.
Let us prove (\ref{TV}). The first inequality of (\ref{TV}) follows
immediately from the definition of the total variation distance. The
second one is a simple consequence
of LeCam's
inequality (\ref{LeCam}).
Indeed, we have, by the triangle inequality,
\begin{eqnarray*}
d_{\mathrm{TV}} \bigl(P_L({\mathbb
X}),P_{\tau}({\mathbb X}) \bigr) \le \sum_{j=0}^{n-1}d_{\mathrm{TV}}
\bigl(P_{\tau_j}({\mathbb X}), P_{\tau
_{j+1}}({\mathbb X}) \bigr),
\end{eqnarray*}
where $\tau_0=L$, $\tau_n=\tau$, and, for $1\le j\le n-1$,
\begin{eqnarray*}
\tau_j=\sum_{t=1}^{j}(X_t-1)_+
\eta_t+\sum_{t=j+1}^{n}(X_t-1)_+
{\mathbb I}_t(w_1).
\end{eqnarray*}
Furthermore, noting that the sums $\tau_j$ and $\tau_{j+1}$ differ only
by one term, we obtain
\begin{eqnarray*}
d_{\mathrm{TV}} \bigl(P_{\tau_j}({\mathbb X}), P_{\tau_{j+1}}({
\mathbb X}) \bigr) 
\le d_{\mathrm{TV}} \bigl(P_{{\mathbb I}_{j+1}}({\mathbb X}),
P_{\eta
_{j+1}}({\mathbb X}) \bigr).
\end{eqnarray*}
Finally, invoking the inequality
$d_{\mathrm{TV}} (P_{{\mathbb I}_{j+1}}({\mathbb X}), P_{\eta
_{j+1}}({\mathbb X}) )
\le
2m^{-2}X_{j+1}^2$ [see (\ref{LeCam})], we arrive at~(\ref{TV}).

The proof of the second bound\vspace*{1pt} $I_2=o(1)$ is based on the fact that $\eta
\to\Lambda$ and $Y_1\to{\tilde Z}_1$ in probability
as $m\to\infty$.
Details of the proof are given in the separate Lemma~\ref{technicallema}.

\textit{Part} 2.
We write, by inclusion--exclusion,
$ L\ge d \ge L-T$,
where the number
\[
T = \sum_{1\le i<j\le n}\sum_{w\in W\setminus\{w_1\}}
{\mathbb I}_i(w_1){\mathbb I}_j(w_1){
\mathbb I}_i(w){\mathbb I}_j(w)
\]
is at least as large as the number of vertices $w$ having two or more
links to $w_1$. Hence, we have
$\PP(d\not=L)\le\PP(T\ge1)$. In order to prove (\ref{d=L}), we show
that $\PP(T\ge1)\le\varepsilon+o(1)$ for every $\varepsilon\in(0,1)$.
For this purpose, we write
\begin{eqnarray*}
\PP(T\ge1) = \E\PP(T\ge1|{\mathbb X})\le\varepsilon+\PP \bigl( \PP(T\ge 1|{
\mathbb X})\ge\varepsilon \bigr) \le \varepsilon+\PP(\E_{\mathbb X}T\ge
\varepsilon)
\end{eqnarray*}
%
and invoke the bound
$\E_{\mathbb X}T=o_P(1)$. Let us prove this bound.
We write
\begin{eqnarray*}
\E_{\mathbb X}T=\sum_{1\le i<j\le n}(n-1)
\frac{(X_i)_2}{(m)_2}\frac
{(X_j)_2}{(m)_2} \le (n-1)m^{-4}
\bigl(X_1^2+\cdots+X_n^2
\bigr)^2
\end{eqnarray*}
and use the bound $n^{-3/2}(X_1^2+\cdots+X_n^2)=o_P(1)$, which follows
from (\ref{ThirdFact})
(applied to random variables $X_i^{4/3}$, $1\le i\le n$, and $\alpha=3/2$).
\end{pf*}

\begin{lem}\label{technicallema} Assume that the conditions of Theorem
\ref{PTh-0} hold. For $\xi$ defined in~(\ref{tau-xi}) and
every real~$t$, we have
$\lim_{m\to+\infty}\E e^{it\xi}= \E e^{itd_*}$.
\end{lem}

\begin{pf} We write $I_2=\E
e^{it\xi}-\E e^{itd_*}$.

The cases $\E Z=0$ and $\E Z>0$ are considered separately. In the case
where $\E Z=0$,
the random variable
$d_*$ is degenerate, $\PP(d_*=0)=1$ and the lemma follows from the
relation 
\begin{eqnarray*}
\PP(\xi=0) \ge \PP(\eta=0) = \E\PP(\eta=0|{\mathbb X})=\E e^{-S(X)/m}=1-o(1).
\end{eqnarray*}
In the last step, we combined the bound $n^{-1}S(X)-\E Z=o_P(1)$ [see
(\ref{SecondFact})] with $\E Z=0$.
%

Next, we consider the case where $\E Z>0$. We need some more notation. Denote
%
\begin{eqnarray*}
 f(t)&=&\E e^{itd_*},\qquad f_{\mathbb X}(t)=\E_{\mathbb X}e^{itY_1},\qquad
f_*(t)=\E e^{it{\tilde Z}_1},\\[3pt]
 \lambda_{\mathbb X}&=&\E_{\mathbb X}\eta,
\lambda_*=\E\Lambda,
\\[3pt]
 \Delta&=&\E_{\mathbb X} e^{it \xi}-f(t),\qquad \delta=
\bigl(f_{\mathbb X}(t)-1\bigr)\lambda_{\mathbb X}-\bigl(f_*(t)-1\bigr)
\lambda_*
\end{eqnarray*}
and write the Fourier transforms of the compound Poisson distributions
$P_{d_*}$ and $P_{\xi}({\mathbb X})$ in the form
$\E e^{itd_*}=e^{(f_*(t)-1)\lambda_*}$ and
$\E_{\mathbb X} e^{it \xi}=e^{(f_{\mathbb X}(t)-1)\lambda_{\mathbb
X}}$, respectively.
%
%
Introduce the events
\begin{eqnarray*}
\mathcal{ A}(\varepsilon)= \bigl\{ \bigl|n^{-1}S(X)-\E Z\bigr|<2^{-1}
\varepsilon\min \{1,\beta, \E Z\} \bigr\},\qquad \varepsilon>0, %
\end{eqnarray*}
and the function $t\to a(t)=\sup_m \{(\E X_1)^{-1}\E X_1{\mathbb
I}_{\{X_1\ge t\}} \}$.
It follows from Lemma~\ref{auxiliarylema} that conditions (vii) and
(viii) imply that $n^{-1}S(X)-\E Z=o_P(1)$ as $m\to\infty$
and $a(t)=o(1)$ as $t\to+\infty$.
Hence, we have
%
\begin{equation}
\label{A-event} \forall\varepsilon>0\qquad \PP \bigl(\mathcal{ A}(\varepsilon)
\bigr)=1-o(1) \qquad{\mbox{as }} m\to\infty.
\end{equation}
In addition, for every $\varepsilon>0$, we can choose a positive
integer $t_{\varepsilon}$ such that
%
\begin{equation}
\label{t-epsilon} a(t_{\varepsilon})<\varepsilon\quad {\mbox{and}}\quad (\E
Z)^{-1}\E Z {\mathbb I}_{\{Z\ge t_{\varepsilon}\}}<\varepsilon.
\end{equation}

In order to prove $I_2=o(1)$, we show that there exists a number $c>0$,
depending only on $\E Z$ and $\beta$, such that, for any $\varepsilon
\in(0,1)$,
%
\begin{equation}
\label{I-2-c} \limsup_n|I_2|\le c\varepsilon.
\end{equation}
In the proof of (\ref{I-2-c}), we assume that $m,n $ are so large that
$\beta\le2m/n\le4\beta$ and $\E X_1\le2\E Z\le4\E X_1$.
In particular, on the event $\mathcal{ A}(\varepsilon)$, we have
%
\begin{equation}
\label{m-is-large} m^{-1}S(X)\le3\beta^{-1}\E Z \quad{\mbox{and}}\quad
S(X)>2^{-1}n\E Z\ge4^{-1}n\E X_1.
\end{equation}
We fix $\varepsilon$ and
write
%
\begin{equation}
\qquad I_2=\E\Delta=I_{21}+I_{22},\qquad I_{21}=
\E\Delta{\mathbb I}_{\mathcal{ A}(\varepsilon)},\qquad I_{22}=\E\Delta(1-{\mathbb
I}_{\mathcal{ A}(\varepsilon)}).
\end{equation}
Note that $|\Delta|\le2$ (as the absolute value of the Fourier
transform of a probability distribution is at most $1$).
This inequality, together with (\ref{A-event}), implies
$I_{22}=o(1)$ as $m\to\infty$. Next, we estimate $I_{21}$. Combining
the identity
$\Delta=f(t) (e^{\delta}-1 )$
with the inequalities $|f(t)|\le1$ and $|e^{\delta}-1|\le\sum_{k\ge
1}|\delta|^k/k!\le|\delta|e^{|\delta|}$, we obtain
\begin{eqnarray*}
|I_{21}|\le\E|\delta|e^{|\delta|} {\mathbb I}_{\mathcal{ A}(\varepsilon)} \le
e^{8\lambda_*}\E|\delta| {\mathbb I}_{\mathcal{ A}(\varepsilon)}.
\end{eqnarray*}
In the last step, we used the bound $|\delta|\le8\lambda_*$, which
follows from
the inequality
$|\delta|\le2\lambda_{\mathbb X}+2\lambda_*$ and the inequality
$\lambda_{\mathbb X}\le3\lambda_*$; see (\ref{m-is-large}).

Next, we show that $\E|\delta| {\mathbb I}_{\mathcal{ A}(\varepsilon)}\le
(2+7\lambda_*)\varepsilon+o(1)$.
We write
\begin{eqnarray*}
\delta=\bigl(f_{\mathbb X}(t)-1\bigr) (\lambda_{\mathbb X}-\lambda_*)+
\bigl(f_{\mathbb
X}(t)-f_*(t)\bigr)\lambda_*,
\end{eqnarray*}
estimate $|\delta|\le2|\lambda_{\mathbb X}-\lambda_*|+\lambda_*|
f_{\mathbb X}(t)-f_*(t)|$ and substitute
\begin{eqnarray*}
|\lambda_{\mathbb X}-\lambda_*|\le\bigl|n^{-1}S(X)-\E
Z\bigr|nm^{-1} + \bigl|nm^{-1}-\beta^{-1}\bigr| \E Z.
\end{eqnarray*}
In this way, we obtain
$\E|\delta| {\mathbb I}_{\mathcal{ A}(\varepsilon)}\le2 I_{31}+2
I_{32}+\lambda_*I_{33}$,
where
%
\begin{eqnarray}
\label{I-123}  I_{31}&=&nm^{-1}\E\bigl\llvert
n^{-1}S(X)-\E Z\bigr\rrvert {\mathbb I}_{\mathcal{
A}(\varepsilon)}\le\varepsilon,
\nonumber\\
I_{32}&=&\bigl|nm^{-1}-\beta^{-1}\bigr| \E
Z=o(1),
\\
\nonumber
 I_{33}&=&\E\bigl\llvert f_{\mathbb X}(t)-f_*(t)\bigr
\rrvert {\mathbb I}_{\mathcal{
A}(\varepsilon)}.
\end{eqnarray}
It remains to estimate $I_{33}$. We expand $f_{\mathbb
X}(t)-f_*(t)=R_1+R_2+R_3$, where
%
\begin{eqnarray}
\nonumber
 R_1&=&\sum_{0\le k\le t_{\varepsilon}}e^{itk}
\delta_P(k), \qquad\delta_P(k)=\PP\bigl(Y_1({
\mathbb X})=k|{\mathbb X}\bigr)-\PP({\tilde Z}_1=k),
\\
\label{R2-suvalgytas} \qquad |R_2|&\le&\PP({\tilde Z}_1\ge
t_{\varepsilon}+1) = (\E Z)^{-1}\E Z{\mathbb I}_{\{Z\ge t_{\varepsilon}+2\}}<
\varepsilon,
\\
\nonumber
|R_3|&\le&\PP \bigl( Y_1({\mathbb X})\ge
t_{\varepsilon}+1|{\mathbb X} \bigr) = \bigl(S(X)\bigr)\sum
_{j=1}^nX_j{\mathbb I}_{\{X_j>t_{\varepsilon}+2\}}
\end{eqnarray}
and observe that the last inequality of (\ref{m-is-large}) implies
%
\begin{equation}
\label{R-3-suvalgytas}  \E|R_3|{\mathbb I}_{\mathcal{ A}(\varepsilon)} \le 4(\E
X_1)^{-1}\E \Biggl(n^{-1}\sum
_{j=1}^nX_j{\mathbb I}_{\{
X_j>t_{\varepsilon}+2\}}
\Biggr) \le 4a(t_{\varepsilon}+2)<4\varepsilon.\hspace*{-35pt}
\end{equation}
Finally, we write $|R_1|\le\sum_{0\le k\le t_{\varepsilon}}|\delta_P(k)|$ and expand $\delta_P(k)=\sum_{j=1}^3r_j(k+1)$, where
\begin{eqnarray*}
r_1(k)=\frac{k}{S(X)} \bigl(n_k-np_x(k)
\bigr),\qquad r_2(k)=kp_x(k) \biggl(\frac{n}{S(X)}-
\frac{1}{\E Z} \biggr),
\end{eqnarray*}
and $r_3(k)=\frac{k}{\E Z} (p_x(k)-p_z(k) )$. Here we denote
$p_x(k)=\PP(X_1=k)$ and $p_z(k)=\PP(Z=k)$.
Since the number $t_{\varepsilon}$ is fixed, it follows from (vii) that
%
\begin{equation}
\label{r-3-suvalgytas} \sum_{0\le k\le t_{\varepsilon}}\bigl|r_3(k+1)\bigr|=o(1)
\qquad{\mbox{as }} m\to\infty.
\end{equation}
Furthermore, on the event $\mathcal{ A}(\varepsilon)$, we obtain from (\ref
{m-is-large}) that
\begin{eqnarray*}
\bigl|r_1(k)\bigr|\le4(\E X_1)^{-1}k\biggl\llvert
\frac{n_k}{n}-p_k(x)\biggr\rrvert ,\qquad \bigl|r_2(k)\bigr|\le2
\varepsilon(\E X_1)^{-1}kp_k(x).
\end{eqnarray*}
Hence, we have $\sum_{0\le k\le t_{\varepsilon}}|r_2(k+1)|{\mathbb
I}_{\mathcal{ A}(\varepsilon)}\le2\varepsilon$ and
\begin{eqnarray*}
\E\sum_{0\le k\le t_{\varepsilon}}\bigl|r_1(k+1)\bigr|{
\mathbb I}_{\mathcal{
A}(\varepsilon)} &\le & \frac{4}{\E X_1}\sum
_{1\le k\le t_{\varepsilon}+1}k \biggl(\E \biggl(\frac
{n_k}{n}-p_k(x)
\biggr)^2 \biggr)^{1/2}
\\
& \le & \frac{4}{\E X_1}n^{-1/2}\sum
_{1\le k\le t_{\varepsilon}+1}k
\\
&=& o(1)
\end{eqnarray*}
as $m\to\infty$. We conclude that $\E|R_1|{\mathbb I}_{\mathcal{
A}(\varepsilon)}\le2\varepsilon+o(1)$.
This bound,
together with~(\ref{R2-suvalgytas}) and (\ref{R-3-suvalgytas}),
shows the bound $I_{33}\le7\varepsilon+o(1)$. The proof is complete.
\end{pf}

\begin{pf*}{Proof of Lemma~\ref{PLemma-0}}
Let us prove (\ref{n*}). We only consider the case where $\varepsilon
m(\max\{1, n_*\})^{-1}\ge2$, since otherwise
(\ref{n*}) is obvious.
Given $0< k\le2^{-1}\varepsilon m$, introduce the random variables
%
\begin{eqnarray*}
d_A(k)=\sum_{i=1}^n{\mathbb
I}_{\{w_1\in D_i\}} {\mathbb I}_{\{2\le
X_i<\varepsilon m/k\}}, \qquad d_B(k)=\sum
_{i=1}^n{\mathbb I}_{\{w_1\in D_i\}} {
\mathbb I}_{\{X_i\ge
\varepsilon m/k\}}
\end{eqnarray*}
and note that the equality of events $\{d\ge1\}=\{d_A(k)\ge1\}\cup\{
d_B(k)\ge1\}$ implies
%
\begin{equation}
\label{A--B} 
\PP\bigl(d_B(k)\ge1\bigr)\ge\PP(d\ge1)-\PP
\bigl(d_A(k)\ge1\bigr).
\end{equation}
%
Combining (\ref{A--B}) with the inequalities $\PP(d\ge\varepsilon
mk^{-1}-1)\ge
\PP(d_B(k)\ge1)$
and
%
\begin{eqnarray*}
\PP\bigl(d_A(k)\ge1\bigr) \le \E \bigl(\E\bigl(d_A(k)|X_1,
\ldots, X_n\bigr) \bigr) = \E\sum_{i=1}^nX_im^{-1}{
\mathbb I}_{\{2\le X_i<\varepsilon m/k\}} < \varepsilon n_* k^{-1},
\end{eqnarray*}
we obtain $\PP(d\ge\varepsilon mk^{-1}-1)\ge\PP(d\ge1)-\varepsilon
n_* k^{-1}$.
Finally, we choose $k=\max\{1,n_*\}$ and obtain (\ref{n*}).\vadjust{\goodbreak}

The proof of (\ref{m>>n}) is much the same, but now we choose $k=n$.

Let us prove (\ref{case-b}). Given
$N=\sum_{i=1}^n{\mathbb I}_{\{X_i\ge2\}}$, let $Y_N$ be a random
variable with binomial distribution
$\Bin(N, 2m^{-1})$. We have
\begin{eqnarray*}
\PP(d\ge C)\ge\PP(Y_N\ge C)\ge\PP\bigl(Y_N\ge C|N
\ge2^{-1}n_*\bigr)\PP\bigl(N\ge 2^{-1}n_*\bigr)=1-o(1).
\end{eqnarray*}
In the last step, we used the simple facts about the binomial distribution
that $\E N=n_*\to\infty$ implies
$\PP(N\ge2^{-1}n_*)=1-o(1)$ and that $\E(Y_N|N=2^{-1}n_*)=
n_*m^{-1}\to+\infty$
implies
\begin{eqnarray*}
\PP\bigl(Y_N\ge C|N\ge2^{-1}n_*\bigr)\ge\PP
\bigl(Y_N\ge C|N= 2^{-1}n_*\bigr)=1-o(1).\hspace*{35pt}\qed
\end{eqnarray*}
\noqed\end{pf*}

\begin{pf*}{Proof of Theorem~\ref{PTh-1}}
The convergence of the degree distribution follows from Theorem \ref
{PTh-0}. Assuming, in addition,
that $\lim_{m\to\infty}\E X^3_{1}=\E Z^3<\infty$, we obtain
the convergence of moments $\E(X_1)_k\to\E(Z)_k$, $k=2,3$, as $m,n\to
\infty$.
Hence, we can write (\ref{alpha-*}) in the form
\[
\alpha^*(n,m,P)=\frac{ \E(Z)_3}{\beta^{-1} (\E(Z)_2 )^2+\E
(Z)_3}+o(1).
\]
Invoking the identities $\E(Z)_2=\beta  \E d_*$ and
$\E(Z)_3=\beta (\E(d_*)_2-(\E d_*)^2 )$, we obtain (\ref{alpha-*+}).

In the case where $\E Z^2<\infty$ and $\E Z^3=\infty$,
we have $\E(X_1)_2\to\E(Z)_2<\infty$ and $\E(X_1)_3\to\infty$.
Therefore, (\ref{alpha-*})
implies $\alpha^*\to1$.
\end{pf*}

\begin{pf*}{Proof of Theorem~\ref{PTh-2}} In the proof, we use some
new notation. Let
$\{w_2^*, w_3^*\}$ be a random pair of vertices uniformly distributed
in the set of all pairs from $W\setminus\{w_1\}$.
Let ${\overline\delta}=(\delta_1,\delta_2,\ldots, \delta_n)$ denote an
ordered collection
of subsets of $W$. We call ${\overline\delta}$
a set-valued vector.
Introduce the random set-valued vector ${\overline D}=(D_1^*,D_2^*,\ldots
,D_n^*)$, where
$D_i^*=D_i$ for $w_1\in D_i$ and $D_i^*=\varnothing$ otherwise.
Denote ${\overline X}=(X_1^*,\ldots, X_n^*)$,
where $X_i^*=|D_i^*|$.
Introduce the function
${\overline\delta}\to h({\overline\delta})
=
2\sum_{i=1}^n\bigl({|\delta_i|-1\atop 2}\bigr){\mathbb I}_{\{\delta_i\not
=\varnothing\}}$
and the random variable
$H=h({\overline D})=\sum_{i=1}^n(X_i^*-1)_2{\mathbb I}_{\{X_i^*\ge1\}}$.
A collection $\mathcal{ I}$ of set-valued vectors is identified with
the event ${\overline D}\in\mathcal{ I}$.
We denote
$\mathcal{ S}_\mathcal{ I}
=
\E h({\overline D}){\mathbb I}_\mathcal{ I}
=
\sum_{{\overline\delta}\in\mathcal{ I}}h(\overline\delta)\PP({\overline
D}={\overline\delta})$.
By $\PP_{\overline\delta}$ we denote the conditional probability given
the event $\{{\overline D}={\overline\delta}\}$.
Denote $n_{\overline\delta}=\sum_{i=1}^n{\mathbb I}_{\{\delta
_i=\varnothing\}}$.

Now we fix an integer $k\ge2$ such that $\PP(d_*=k)>0$ and assume that
$m,n$ are so large
that $\PP(d=k)\ge2^{-1}\PP(d_*=k)>0$.
Introduce the events
\begin{eqnarray*}
\Delta=\bigl\{w_1\sim w^*_2, w_1\sim
w^*_3, w^*_2\sim w^*_3\bigr\},\qquad\hspace*{-1pt} \mathcal{ A}=
\{d=k\},\qquad\hspace*{-1pt} \mathcal{ A}_k=\{d=L, L=k\} 
\end{eqnarray*}
and the probabilities
\begin{eqnarray*}
p_{\star}&=&\PP\bigl(w_1\sim w^*_2,
w_1\sim w^*_3, d=k\bigr),\\
 {\tilde p}_{\star}&=&\PP
\bigl(w_1\sim w^*_2, w_1\sim
w^*_3, d=k, d\not=L\bigr).
\end{eqnarray*}
Here $\{w_2^*, w_3^*\}$ is a random subset of $W\setminus\{w_1\}$
uniformly distributed in the class of subsets of size $2$. We assume
that $\{w_2^*, w_3^*\}$
is independent
of the random sets $D_1,\ldots, D_n$ defining the intersection graph.

In order to prove (\ref{alpha-*++}), we write $\alpha^{*[k]}$ in the form
\[
\alpha^{*[k]}=\frac{\PP(\Delta\cap\mathcal{ A})}{p_{\star}}
\]
and invoke the identities
%
\begin{eqnarray}
\label{D-A-22}  \PP(\Delta\cap\mathcal{ A})&=&\frac{1}{(m-1)_2}\E
d_{2*}{\mathbb I}_{\{
d_*=k\}}+o\bigl(m^{-2}\bigr),
\\
\label{p-star} p_{\star}&=&{\pmatrix{k
\cr
2}} {\pmatrix{m-1
\cr
2}}^{-1}\PP(d_*=k)+o\bigl(m^{-2}\bigr).
\end{eqnarray}
In order to show (\ref{p-star}), we note that
every pair $\{w',w''\}\subset W\setminus\{w_1\}$ has the same probability
to be covered by the neighborhood of $w_1$. Hence, by symmetry we have
%
\begin{eqnarray}
\label{p-p-star} p_{\star}&=&{\pmatrix{k
\cr
2}} {\pmatrix{m-1
\cr
2}}^{-1}\PP(d=k),
\nonumber
\\[-8pt]
\\[-8pt]
\nonumber
  {\tilde p}_{\star}&=&{\pmatrix{k
\cr
2}} {
\pmatrix{m-1
\cr
2}}^{-1}\PP(d=k, d\not=L).
\end{eqnarray}
Now (\ref{p-star}) follows from the first identity and the convergence
$\PP(d=k)\to\PP(d_*=k)$.

Let us show (\ref{D-A-22}). Observe that
$\PP(\Delta\cap\mathcal{ A}\cap\{d\not= L\})\le{\tilde p}_{\star}$, and by
(\ref{p-p-star}) and~(\ref{d=L}) we have ${\tilde p}_{\star
}=o(m^{-2})$. Hence,
%
\begin{equation}
\label{DD+} \PP(\Delta\cap\mathcal{ A})=\PP(\Delta\cap\mathcal{
A}_k)+o\bigl(m^{-2}\bigr).
\end{equation}

Next, we expand, by the total probability formula,
%
\begin{equation}
\label{DDDAB} \PP(\Delta\cap\mathcal{ A}_k) = \sum
_{{\overline\delta}\in\mathcal{ A}_k} \PP_{\overline\delta}(\Delta)\PP({\overline D}={\overline
\delta})
\end{equation}
and note that
$\PP_{\overline\delta}(\mathcal{ D}_{\overline\delta})
\le
\PP_{\overline\delta}(\Delta)
\le
\PP_{\overline\delta}(\mathcal{ C}_{\overline\delta})+\PP_{\overline
\delta}(\mathcal{ D}_{\overline\delta})$,
where
\begin{eqnarray*}
\mathcal{ C}_{\overline\delta} = \bigl\{ \exists D_j\dvtx
w^*_2, w^*_3\in D_j {\mbox{ and }}
w_1\notin D_j\bigr\}, \qquad \mathcal{ D}_{\overline\delta} =
\bigl\{\exists D_j\dvtx w_1,w^*_2,w^*_3
\in D_j\bigr\}.
\end{eqnarray*}
Observe that
%
\begin{equation}
\label{DDDL} \PP_{\overline\delta}(\mathcal{ D}_{\overline\delta})=\frac{h({\overline
\delta})}{(m-1)_2}
\end{equation}
and, since $n_{\overline\delta}<n$,
%
\begin{equation}
\label{DDDC} \PP_{\overline\delta}(\mathcal{ C}_{\overline\delta}) = 1- \biggl( 1-
\frac{\E(X')_2}{(m-1)_2} \biggr)^{n_{\overline\delta}} < n\frac{\E(X')_2}{(m-1)_2}=O
\bigl(m^{-1}\bigr).
\end{equation}
Here $X'$ denotes the random variable $|D_i|$ conditioned on the event
$w_1\notin D_i$.
Collecting (\ref{DDDL}) and (\ref{DDDC}) in (\ref{DDDAB}), we obtain
%
\begin{equation}
\label{DELTA-A-B} \PP(\Delta\cap\mathcal{ A}_k) = \frac{1}{(m-1)_2}
\mathcal{ S}_{\mathcal{ A}_k}+O\bigl(m^{-3}\bigr).
\end{equation}
Now we replace the event $\mathcal{ A}_k$ by $\mathcal{ A}$ in
(\ref{DELTA-A-B}). For the left-hand side, we apply~(\ref{DD+}). For
the right-hand side, we
apply the inequalities
%
\begin{equation}
\label{SAB} \mathcal{ S}_{\mathcal{ A}}\ge\mathcal{ S}_{\mathcal{ A}_k} \ge
\mathcal{ S}_{\mathcal{ A}}-\pmatrix{k
\cr
2}\PP(d=k, d\not=L) = \mathcal{
S}_{\mathcal{ A}}-o(1)
\end{equation}
[here we used $h({\overline\delta})\le\bigl({k\atop  2}\bigr)$ and (\ref{d=L})].
We obtain
\[
\PP(\Delta\cap\mathcal{ A})=\frac{1}{(m-1)_2}\mathcal{ S}_{\mathcal{ A}}+o
\bigl(m^{-2}\bigr).
\]
Finally, (\ref{D-A-22}) follows by the convergence of $\mathcal{ S}_{\mathcal{
A}}=\E H{\mathbb I}_{\{d=k\}}$ to
$\E d_{2*}{\mathbb I}_{\{d_*=k\}}$
as $m,n\to\infty$. We derive this convergence from the weak convergence
of bivariate random vectors $(H,L)\to(d_{2*}, d_*)$, which is
obtained using the same argument as that of the proof of
Theorem~\ref{PTh-0}.
\end{pf*}

\section*{Acknowledgments}
I would like to thank anonymous referees for their valuable comments
and suggestions.



\printaddresses


\begin{thebibliography}{32}

\bibitem{Barbour2011}
\begin{barticle}[mr]
\bauthor{\bsnm{Barbour},~\bfnm{A.~D.}\binits{A.~D.}} \AND
\bauthor{\bsnm{Reinert},~\bfnm{G.}\binits{G.}}
(\byear{2011}).
\btitle{The shortest distance in random multi-type intersection graphs}.
\bjournal{Random Structures Algorithms}
\bvolume{39}
\bpages{179--209}.
\bid{doi={10.1002/rsa.20351}, issn={1042-9832}, mr={2850268}}
\bptok{imsref}%
\end{barticle}
\endbibitem

\bibitem{Barrat2000}
\begin{barticle}[auto:STB|2012/11/23|13:23:43]
\bauthor{\bsnm{Barrat},~\bfnm{A.}\binits{A.}} \AND
\bauthor{\bsnm{Weigt},~\bfnm{M.}\binits{M.}}
(\byear{2000}).
\btitle{On the properties of small-world networks}.
\bjournal{Eur. Phys. J. B}
\bvolume{13}
\bpages{547--560}.
\bptok{imsref}%
\end{barticle}
\endbibitem

\bibitem{Blackburn2009}
\begin{barticle}[mr]
\bauthor{\bsnm{Blackburn},~\bfnm{Simon~R.}\binits{S.~R.}} \AND
\bauthor{\bsnm{Gerke},~\bfnm{Stefanie}\binits{S.}}
(\byear{2009}).
\btitle{Connectivity of the uniform random intersection graph}.
\bjournal{Discrete Math.}
\bvolume{309}
\bpages{5130--5140}.
\bid{doi={10.1016/j.disc.2009.03.042}, issn={0012-365X}, mr={2548914}}
\bptok{imsref}%
\end{barticle}
\endbibitem

\bibitem{Bloznelis2008}
\begin{barticle}[mr]
\bauthor{\bsnm{Bloznelis},~\bfnm{M.}\binits{M.}}
(\byear{2008}).
\btitle{Degree distribution of a typical vertex in a general random
intersection graph}.
\bjournal{Lith. Math. J.}
\bvolume{48}
\bpages{38--45}.
\bid{doi={10.1007/s10986-008-0004-7}, issn={0363-1672}, mr={2398169}}
\bptok{imsref}%
\end{barticle}
\endbibitem

\bibitem{Bloznelis2010}
\begin{barticle}[mr]
\bauthor{\bsnm{Bloznelis},~\bfnm{Mindaugas}\binits{M.}}
(\byear{2010}).
\btitle{A random intersection digraph: Indegree and outdegree distributions}.
\bjournal{Discrete Math.}
\bvolume{310}
\bpages{2560--2566}.
\bid{doi={10.1016/j.disc.2010.06.018}, issn={0012-365X}, mr={2669380}}
\bptok{imsref}%
\end{barticle}
\endbibitem

\bibitem{Bloznelis2010a}
\begin{barticle}[mr]
\bauthor{\bsnm{Bloznelis},~\bfnm{Mindaugas}\binits{M.}}
(\byear{2010}).
\btitle{The largest component in an inhomogeneous random intersection graph
with clustering}.
\bjournal{Electron. J. Combin.}
\bvolume{17}
\bpages{Research Paper 110, 17}.
\bid{issn={1077-8926}, mr={2679564}}
\bptok{imsref}%
\end{barticle}
\endbibitem

\bibitem{Breiman}
\begin{bbook}[mr]
\bauthor{\bsnm{Breiman},~\bfnm{Leo}\binits{L.}}
(\byear{1968}).
\btitle{Probability}.
\bpublisher{Addison-Wesley}, \blocation{Reading, MA.}
\bid{mr={0229267}}
\bptok{imsref}%
\end{bbook}
\endbibitem

\bibitem{Britton2008}
\begin{barticle}[mr]
\bauthor{\bsnm{Britton},~\bfnm{Tom}\binits{T.}},
\bauthor{\bsnm{Deijfen},~\bfnm{Maria}\binits{M.}},
\bauthor{\bsnm{Lager{\aa}s},~\bfnm{Andreas~N.}\binits{A.~N.}} \AND
\bauthor{\bsnm{Lindholm},~\bfnm{Mathias}\binits{M.}}
(\byear{2008}).
\btitle{Epidemics on random graphs with tunable clustering}.
\bjournal{J. Appl. Probab.}
\bvolume{45}
\bpages{743--756}.
\bid{doi={10.1239/jap/1222441827}, issn={0021-9002}, mr={2455182}}
\bptok{imsref}%
\end{barticle}
\endbibitem

\bibitem{Deijfen}
\begin{barticle}[mr]
\bauthor{\bsnm{Deijfen},~\bfnm{Maria}\binits{M.}} \AND
\bauthor{\bsnm{Kets},~\bfnm{Willemien}\binits{W.}}
(\byear{2009}).
\btitle{Random intersection graphs with tunable degree distribution and
clustering}.
\bjournal{Probab. Engrg. Inform. Sci.}
\bvolume{23}
\bpages{661--674}.
\bid{doi={10.1017/S0269964809990064}, issn={0269-9648}, mr={2535025}}
\bptok{imsref}%
\end{barticle}
\endbibitem

\bibitem{eschenauer2002}
\begin{bincollection}[auto:STB|2012/11/23|13:23:43]
\bauthor{\bsnm{Eschenauer},~\bfnm{L.}\binits{L.}} \AND
\bauthor{\bsnm{Gligor},~\bfnm{V.~D.}\binits{V.~D.}}
(\byear{2002}).
\btitle{A key-management scheme for distributed sensor networks}.
In \bbooktitle{Proceedings of the 9th ACM Conference on Computer and Communications Security, Washington, DC}
\bpages{41--47}.
\bptok{imsref}%
\end{bincollection}
\endbibitem


\bibitem{Foss}
\begin{bbook}[mr]
\bauthor{\bsnm{Foss},~\bfnm{Sergey}\binits{S.}},
\bauthor{\bsnm{Korshunov},~\bfnm{Dmitry}\binits{D.}} \AND
\bauthor{\bsnm{Zachary},~\bfnm{Stan}\binits{S.}}
(\byear{2011}).
\btitle{An Introduction to Heavy-Tailed and Subexponential Distributions}.
\bpublisher{ACM}, \blocation{New York}.
\bid{doi={10.1007/978-1-4419-9473-8}, mr={2810144}}
\bptnote{check year}%
\bptok{imsref}%
\end{bbook}
\endbibitem

\bibitem{Foudalis2011}
\begin{bincollection}[mr]
\bauthor{\bsnm{Foudalis},~\bfnm{Ilias}\binits{I.}},
\bauthor{\bsnm{Jain},~\bfnm{Kamal}\binits{K.}},
\bauthor{\bsnm{Papadimitriou},~\bfnm{Christos}\binits{C.}} \AND
\bauthor{\bsnm{Sideri},~\bfnm{Martha}\binits{M.}}
(\byear{2011}).
\btitle{Modeling social networks through user background and behavior}.
In \bbooktitle{Algorithms and Models for the Web Graph}.
\bseries{Lecture Notes in Computer Science}
\bvolume{6732}
\bpages{85--102}.
\bpublisher{Springer}, \blocation{Heidelberg}.
\bid{doi={10.1007/978-3-642-21286-4_8}, mr={2842315}}
\bptok{imsref}%
\end{bincollection}
\endbibitem

\bibitem{godehardt2001}
\begin{barticle}[mr]
\bauthor{\bsnm{Godehardt},~\bfnm{Erhard}\binits{E.}} \AND
\bauthor{\bsnm{Jaworski},~\bfnm{Jerzy}\binits{J.}}
(\byear{2001}).
\btitle{Two models of random intersection graphs and their applications}.
\bjournal{Electron. Notes Discrete Math.}
\bvolume{10}
\bpages{129--132}.
\bptok{imsref}%
\end{barticle}
\endbibitem

\bibitem{godehardt2003}
\begin{bincollection}[mr]
\bauthor{\bsnm{Godehardt},~\bfnm{E.}\binits{E.}} \AND
\bauthor{\bsnm{Jaworski},~\bfnm{J.}\binits{J.}}
(\byear{2003}).
\btitle{Two models of random intersection graphs for classification}.
In \bbooktitle{Exploratory Data Analysis in Empirical Research}
\bpages{67--81}.
\bpublisher{Springer}, \blocation{Berlin}.
\bid{mr={2074223}}
\bptok{imsref}%
\end{bincollection}
\endbibitem

\bibitem{GJR2010}
\begin{bincollection}[auto:STB|2012/11/23|13:23:43]
\bauthor{\bsnm{Godehardt},~\bfnm{E.}\binits{E.}},
\bauthor{\bsnm{Jaworski},~\bfnm{J.}\binits{J.}} \AND
\bauthor{\bsnm{Rybarczyk},~\bfnm{K.}\binits{K.}}
(\byear{2012}).
\btitle{Clustering coefficients of random intersection graphs}.
In \bbooktitle{Studies in Classification, Data Analysis and Knowledge Organization}
\bpages{243--253}.
\bpublisher{Springer}, \blocation{Berlin}.
\bptok{imsref}%
\end{bincollection}
\endbibitem

\bibitem{Guillaume+L2004}
\begin{barticle}[mr]
\bauthor{\bsnm{Guillaume},~\bfnm{Jean-Loup}\binits{J.-L.}} \AND
\bauthor{\bsnm{Latapy},~\bfnm{Matthieu}\binits{M.}}
(\byear{2004}).
\btitle{Bipartite structure of all complex networks}.
\bjournal{Inform. Process. Lett.}
\bvolume{90}
\bpages{215--221}.
\bid{doi={10.1016/j.ipl.2004.03.007}, issn={0020-0190}, mr={2054656}}
\bptok{imsref}%
\end{barticle}
\endbibitem

\bibitem{JKS}
\begin{barticle}[mr]
\bauthor{\bsnm{Jaworski},~\bfnm{Jerzy}\binits{J.}},
\bauthor{\bsnm{Karo{\'n}ski},~\bfnm{Micha{\l}}\binits{M.}} \AND
\bauthor{\bsnm{Stark},~\bfnm{Dudley}\binits{D.}}
(\byear{2006}).
\btitle{The degree of a typical vertex in generalized random intersection graph
models}.
\bjournal{Discrete Math.}
\bvolume{306}
\bpages{2152--2165}.
\bid{doi={10.1016/j.disc.2006.05.013}, issn={0012-365X}, mr={2255609}}
\bptok{imsref}%
\end{barticle}
\endbibitem

\bibitem{JaworskiStark2008}
\begin{barticle}[mr]
\bauthor{\bsnm{Jaworski},~\bfnm{Jerzy}\binits{J.}} \AND
\bauthor{\bsnm{Stark},~\bfnm{Dudley}\binits{D.}}
(\byear{2008}).
\btitle{The vertex degree distribution of passive random intersection graph
models}.
\bjournal{Combin. Probab. Comput.}
\bvolume{17}
\bpages{549--558}.
\bid{doi={10.1017/S0963548308009103}, issn={0963-5483}, mr={2433940}}
\bptok{imsref}%
\end{barticle}
\endbibitem

\bibitem{karonski1999}
\begin{barticle}[mr]
\bauthor{\bsnm{Karo{\'n}ski},~\bfnm{Micha{\l}}\binits{M.}},
\bauthor{\bsnm{Scheinerman},~\bfnm{Edward~R.}\binits{E.~R.}} \AND
\bauthor{\bsnm{Singer-Cohen},~\bfnm{Karen~B.}\binits{K.~B.}}
(\byear{1999}).
\btitle{On random intersection graphs: The subgraph problem}.
\bjournal{Combin. Probab. Comput.}
\bvolume{8}
\bpages{131--159}.
\bid{doi={10.1017/S0963548398003459}, issn={0963-5483}, mr={1684626}}
\bptok{imsref}%
\end{barticle}
\endbibitem

\bibitem{Newman2003}
\begin{barticle}[auto:STB|2012/11/23|13:23:43]
\bauthor{\bsnm{Newman},~\bfnm{M.~E.~J.}\binits{M.~E.~J.}}
(\byear{2003}).
\btitle{Properties of highly clustered networks}.
\bjournal{Phys. Rev. E}
\bvolume{68}
\bpages{026121}.
\bptok{imsref}%
\end{barticle}
\endbibitem

\bibitem{Newman2001}
\begin{barticle}[auto:STB|2012/11/23|13:23:43]
\bauthor{\bsnm{Newman},~\bfnm{M.~E.~J.}\binits{M.~E.~J.}},
\bauthor{\bsnm{Strogatz},~\bfnm{S.~H.}\binits{S.~H.}} \AND
\bauthor{\bsnm{Watts},~\bfnm{D.~J.}\binits{D.~J.}}
(\byear{2001}).
\btitle{Random graphs with arbitrary degree distributions and their
applications}.
\bjournal{Phys. Rev. E}
\bvolume{64}
\bpages{026118}.
\bptok{imsref}%
\end{barticle}
\endbibitem

\bibitem{Newman+W+S2002}
\begin{barticle}[auto:STB|2012/11/23|13:23:43]
\bauthor{\bsnm{Newman},~\bfnm{M.~E.~J.}\binits{M.~E.~J.}},
\bauthor{\bsnm{Watts},~\bfnm{D.~J.}\binits{D.~J.}} \AND
\bauthor{\bsnm{Strogatz},~\bfnm{S.~H.}\binits{S.~H.}}
(\byear{2002}).
\btitle{Random graph models of social networks}.
\bjournal{Proc. Natl. Acad. Sci. USA}
\bvolume{99} \bnote{(Suppl. 1)}
\bpages{2566--2572}.
\bptok{imsref}%
\end{barticle}
\endbibitem

\bibitem{Spirakis2011}
\begin{barticle}[mr]
\bauthor{\bsnm{Nikoletseas},~\bfnm{S.}\binits{S.}},
\bauthor{\bsnm{Raptopoulos},~\bfnm{C.}\binits{C.}} \AND
\bauthor{\bsnm{Spirakis},~\bfnm{P.~G.}\binits{P.~G.}}
(\byear{2011}).
\btitle{On the independence number and {H}amiltonicity of uniform random
intersection graphs}.
\bjournal{Theoret. Comput. Sci.}
\bvolume{412}
\bpages{6750--6760}.
\bid{doi={10.1016/j.tcs.2011.09.003}, issn={0304-3975}, mr={2885093}}
\bptok{imsref}%
\end{barticle}
\endbibitem

\bibitem{RBarabasi2003}
\begin{barticle}[auto:STB|2012/11/23|13:23:43]
\bauthor{\bsnm{Ravasz},~\bfnm{L.}\binits{L.}} \AND
\bauthor{\bsnm{Barab{\'a}si},~\bfnm{A.~L.}\binits{A.~L.}}
(\byear{2003}).
\btitle{Hierarchical organization in complex networks}.
\bjournal{Phys. Rev. E}
\bvolume{67}
\bpages{026112}.
\bptok{imsref}%
\end{barticle}
\endbibitem

\bibitem{Rybarczyk2011}
\begin{barticle}[mr]
\bauthor{\bsnm{Rybarczyk},~\bfnm{Katarzyna}\binits{K.}}
(\byear{2011}).
\btitle{Diameter, connectivity, and phase transition of the uniform random
intersection graph}.
\bjournal{Discrete Math.}
\bvolume{311}
\bpages{1998--2019}.
\bid{doi={10.1016/j.disc.2011.05.029}, issn={0012-365X}, mr={2818878}}
\bptok{imsref}%
\end{barticle}
\endbibitem

\bibitem{Rybarczyk-degree2011}
\begin{bincollection}[auto]
\bauthor{\bsnm{Rybarczyk},~\bfnm{Katarzyna}\binits{K.}}
(\byear{2012}).
\btitle{The degree distribution in random intersection graphs}.
In \bbooktitle{Challenges at the Interface of Data Analysis, Computer Science, and Optimization}
\bpages{291--299}.
\bpublisher{Springer}, \blocation{Berlin}.
\bptok{imsref}%
\end{bincollection}
\endbibitem

\bibitem{stark2004}
\begin{barticle}[mr]
\bauthor{\bsnm{Stark},~\bfnm{Dudley}\binits{D.}}
(\byear{2004}).
\btitle{The vertex degree distribution of random intersection graphs}.
\bjournal{Random Structures Algorithms}
\bvolume{24}
\bpages{249--258}.
\bid{doi={10.1002/rsa.20005}, issn={1042-9832}, mr={2068868}}
\bptok{imsref}%
\end{barticle}
\endbibitem

\bibitem{Steele}
\begin{barticle}[mr]
\bauthor{\bsnm{Steele},~\bfnm{J.~Michael}\binits{J.~M.}}
(\byear{1994}).
\btitle{Le {C}am's inequality and {P}oisson approximations}.
\bjournal{Amer. Math. Monthly}
\bvolume{101}
\bpages{48--54}.
\bid{doi={10.2307/2325124}, issn={0002-9890}, mr={1252705}}
\bptok{imsref}%
\end{barticle}
\endbibitem

\bibitem{storgatz1998}
\begin{barticle}[pbm]
\bauthor{\bsnm{Watts},~\bfnm{D.~J.}\binits{D.~J.}} \AND
\bauthor{\bsnm{Strogatz},~\bfnm{S.~H.}\binits{S.~H.}}
(\byear{1998}).
\btitle{Collective dynamics of ``small-world'' networks}.
\bjournal{Nature}
\bvolume{393}
\bpages{440--442}.
\bid{doi={10.1038/30918}, issn={0028-0836}, pmid={9623998}}
\bptok{imsref}%
\end{barticle}
\endbibitem

\bibitem{Yagan2009}
\begin{bincollection}[auto:STB|2012/11/23|13:23:43]
\bauthor{\bsnm{Yagan},~\bfnm{O.}\binits{O.}} \AND
\bauthor{\bsnm{Makowski},~\bfnm{A.~M.}\binits{A.~M.}}
(\byear{2009}).
\btitle{Random key graphs---Can they be small worlds?}
In \bbooktitle{NETCOM'09: Proceedings of the 2009 First International Conference on Networks \& Communications}
\bpages{313--318}.
\bpublisher{IEEE Computer Society}, \blocation{Washington, DC}.
\bptok{imsref}%
\end{bincollection}
\endbibitem

\end{thebibliography}
\end{document}